\documentclass[pdflatex,sn-mathphys-num]{sn-jnl}%

\usepackage{graphicx}%
\usepackage{multirow}%
\usepackage{amsmath,amssymb,amsfonts}%
\usepackage{amsthm}%
\usepackage{mathrsfs}%
\usepackage{nicefrac}%
\usepackage[title]{appendix}%
\usepackage{xcolor}%
\usepackage{textcomp}%
\usepackage{manyfoot}%
\usepackage{booktabs}%
\usepackage{algorithm, algpseudocode, algorithmicx}
\usepackage{subcaption}%
\usepackage{listings}%
\usepackage{matlab-prettifier}%
\definecolor{MATLABCellColour}{RGB}{252,251,220}%
\usepackage{tikz}
\usetikzlibrary{arrows.meta}
\usetikzlibrary{backgrounds}
\usepackage{pgfplots}
\pgfplotsset{compat=1.18} 
\usepackage{booktabs} 
\usepackage{siunitx}
\usepackage{mathabx}
\definecolor{mycolor1}{rgb}{0.00000,0.44700,0.74100}%
\definecolor{mycolor2}{rgb}{0.85000,0.32500,0.09800}%
\definecolor{mycolor3}{rgb}{0.92900,0.69400,0.12500}%
\usepackage{microtype}
\usepackage{comment}
\theoremstyle{thmstyleone}%
\newtheorem{theorem}{Theorem}%
\newtheorem{proposition}[theorem]{Proposition}%
\newtheorem{lemma}[theorem]{Lemma}%
\newtheorem{corollary}[theorem]{Corollary}%

\theoremstyle{thmstyletwo}%
\newtheorem{example}{Example}%
\newtheorem{remark}{Remark}%

\theoremstyle{thmstylethree}%
\newtheorem{definition}{Definition}%

\begin{document}

\title[Kemeny's constant minimization]{Kemeny's constant minimization for reversible Markov chains via structure-preserving perturbations}

\author[1]{\fnm{Fabio} \sur{Durastante}}\email{fabio.durastante@unipi.it}
\equalcont{These authors contributed equally to this work.}
\author*[1]{\fnm{Miryam} \sur{Gnazzo}}\email{miryam.gnazzo@dm.unipi.it}
\author[1]{\fnm{Beatrice} \sur{Meini}}\email{beatrice.meini@unipi.it}
\equalcont{These authors contributed equally to this work.}
\affil[1]{\orgdiv{Mathematics Department}, \orgname{University of Pisa},
\city{Pisa}, \postcode{56127}, \state{PI}, \country{Italy}}

\abstract{Kemeny’s constant measures the efficiency of a Markov chain in traversing its states. We investigate whether structure-preserving perturbations to the transition probabilities of a reversible Markov chain can improve its connectivity while maintaining a fixed stationary distribution. Although the minimum achievable value for Kemeny’s constant can be estimated, the required perturbations may be infeasible. We reformulate the problem as an optimization task, focusing on solution existence and efficient algorithms, with an {emphasis on} the problem of minimizing Kemeny’s constant under sparsity constraints.}



\maketitle

\section{Introduction}

Kemeny’s constant~\cite{MR636037} quantifies the expected number of steps required for a Markov chain to transition from a starting state to a randomly selected destination state, where the destination is sampled according to the stationary distribution of the Markov chain. This constant serves as a comprehensive measure of the efficiency with which the Markov chain explores its state space. Lower values of Kemeny’s constant indicate faster average transitions between states, enhancing the chain's connectivity, while higher values denote a reduction in communication speed. Kemeny's constant has found applications in diverse areas such as assessing the vulnerability of complex networks via the random walks on their underlying graphs~\cite{LiLinJieShi} and computing edge centrality scores for network analysis~\cite{MR4590023}, developing opinion dynamics models of the evolution of credit scores of enterprises in a social network~\cite{LiangXuChiclana}, and describing epidemic dynamics on networks~\cite{FengHuangBing}. Recently, it has also been utilized to identify graph partitionings that preserve the original graph's time scales by maximizing Kemeny’s constant within a specific partition structure~\cite{doi:10.1021/acs.jpcb.3c08213}.

This work addresses the question of whether it is possible to construct structured perturbations of a Markov chain, specifically targeted modifications to its transition probabilities, that enhance its connectivity by reducing Kemeny's constant. {In this context, it is useful to consider perturbations which preserve the stationary distribution of the starting chain, i.e., preserving the long-term equilibrium behavior of the system, even after modifying the transition probabilities. Similarly, prescribing the structure, i.e., the nonzero pattern of the perturbed transition matrix, reflects practical constraints, since only a limited number of transitions or physical connections can actually be adjusted. Furthermore, requiring a perturbation of \emph{small norm} can also be motivated by implementation costs, since the norm may be interpreted as a measure of the effort, energy expenditure, or resources needed to alter the original dynamics. Thus, our objective is to obtain meaningful reductions in Kemeny's constant through perturbations that are both structurally feasible and quantitatively limited.} 

Optimizing a global measure of Markov chains is closely related to optimizing global measures of complex networks. Examples include the total communicability~\cite{BenziKlymko}, studied in~\cite{MR3439773,MR3479697}; the reduction of a graph's controversy score~\cite{Garimella}; 
the Estrada index~\cite{ESTRADA2000713}---the trace of the matrix exponential of the adjacency matrix---analyzed in~\cite{MasseiTudisco}; and the related natural connectivity investigated in~\cite{MIOBI}. {Following this path, recently in~\cite{Guglielmi} the authors assess the robustness associated with eigenvector centrality via small modifications of complex networks, in \cite{MR4973135} through properties of the Fiedler and Perron eigenpairs, while target values of Katz and PageRank centrality measures are reached via an {Interior Point Method (IPM)} approach in~\cite{Cipolla}; moreover, in~\cite{MR4822696} the authors face the problem of assigning a probability distribution to a given Markov chain.}
{In this research framework, we also treat an application where the Markov chain arises as a random walk on a complex network, with particular emphasis on networks representing national-scale power grids~\cite{MEDJROUBI201714}; see also the analogous application in~\cite{MasseiTudisco}.}

For irreducible Markov chains with a fixed stationary distribution, estimations and bounds for the minimal Kemeny's constant exist \cite{MR3218779,BINI2024137,kirkland2023edgeadditionchangekemenys}. {However, achieving this minimum typically requires modifications of the transition matrix that are not necessarily small in magnitude. In particular, the matrices attaining the lower bound may lie far from the original chain in terms of natural matrix norms, such as the Frobenius norm, meaning that the corresponding perturbations can be too large to be regarded as feasible in applications where only limited adjustments to the dynamics are allowed or where the modification cost must remain controlled.}
The main contribution of this paper is the construction of an algorithmic framework for efficiently determining small-norm perturbations that reduce Kemeny's constant of a given Markov chain, constrained by a given stationary probability vector and predefined sparsity pattern~\cite{Bullo,BulloReview}. %
{Nevertheless, the lack of {smoothness} and convexity of Kemeny's constant over the set of stochastic matrices with given stationary probability vector makes the problem challenging~\cite{BulloReview}. To address the posed question, we reformulate the problem as an optimization task over a specific subset of Markov chains, namely the reversible ones with fixed stationary distribution~\cite{nearestreversiblechain,Reversible}. We then explore the existence of solutions of this optimization problem from a theoretical point of view, and develop efficient optimization algorithms, based on nonlinear programming via IPMs~\cite{MR2115066,MR4865731}, %
and on Riemannian optimization. For the latter, we build a first-order Riemannian geometry over the manifold of stochastic reversible matrices, with fixed stationary probability vector and given nonzero pattern, which is not a straightforward generalization of the tools previously introduced in~\cite{Douik,Reversible}. Both proposed algorithms have been implemented in MATLAB and tested in several numerical examples, and are available on GitHub~\href{https://www.github.com/Cirdans-Home/optimize-kemeny}{Cirdans-Home/optimize-kemeny}.

\subsection{Notation and background material}\label{sec:notation}

To set the stage for the following developments, we first introduce the notation and review fundamental concepts related to discrete-state Markov chains and stochastic matrices. For a more exhaustive exposition, we refer to~\cite{MR410929} and~\cite[\S 8]{Meyer}.
We will denote by $\mathbf{1}$ the vector of all ones, by $D_\mathbf{v}$ the diagonal matrix with diagonal entries the entries of the vector $\mathbf{v}$, by $\operatorname{tr}(A)$ the trace of the matrix $A$, and by $\|A\|_F = \operatorname{tr}(A^\top A)^{\nicefrac{1}{2}}$ the Frobenius norm.

A \textit{nonnegative matrix} $P$ is a real matrix with nonnegative entries and we denote it by $P\geq 0$.
A \textit{stochastic matrix} is a nonnegative matrix $ P $ such that $P \mathbf{1}=\mathbf{1}$. A stochastic matrix has spectral radius 1, and 1 is an eigenvalue. A \textit{probability vector} is a nonnegative vector $\mathbf{v}$ such that $\mathbf{v}^\top \mathbf{1}=1$.

\begin{definition}\label{def:random-walk}
A  \textit{Markov chain} on the space state $\{1,\ldots,n\}$ is a stochastic process $ \{X_h\}_{h \geq 0} $ that satisfies the \textit{Markov property}, 
namely
\begin{align*}
\mathbb{P}(X_{h+1} = i_{h+1} \mid X_h = i_h, X_{h-1} = i_{h-1}, \ldots, X_0 = i_0)  = \mathbb{P}(X_{h+1} = i_{h+1} \mid X_h = i_h),
\end{align*}
for all $ h $ and states $i_j\in \{1,\ldots,n\}$, $j=0,\ldots,h+1$, i.e.,
the future state of the process depends only on the present state and not on the sequence of past states. 
A \textit{homogeneous Markov chain} is also characterized by its \textit{transition probability matrix} $ P = [P_{ij}] {\in \mathbb{R}^{n\times n}} $, where $ P_{ij} = \mathbb{P}(X_{h+1} = j \mid X_h = i) $ represents the probability of transitioning from state $i$ to state $j$, which is a stochastic matrix.
A \textit{stationary distribution}, or \textit{stationary probability vector}, for the Markov chain with transition matrix $P$ is any probability vector $\boldsymbol{\pi}^\top = [\pi_1, \pi_2, \dots, \pi_n]$ such that $\boldsymbol{\pi}^\top P=\boldsymbol{\pi}^\top$. A Markov chain is said to be \textit{irreducible} if its transition matrix $P$ is irreducible{, meaning that for every pair of indices $i$ and $j$ there exists an integer $k > 0$ such that $(P^k)_{ij} > 0$.}
\end{definition}

 An irreducible Markov chain has a unique stationary distribution $\boldsymbol{\pi}$, moreover $\pi_i>0$ for $i=1,\ldots,n$.

\begin{definition}\label{def:reversible-markov-chain}%
An irreducible Markov chain with transition matrix $P$ is said to be \textit{reversible} if the stationary distribution $\boldsymbol{\pi}$ is such that the detailed balance condition holds:
\begin{equation}\label{eq:detailed_balance}
    \pi_i P_{ij} = \pi_j P_{ji}, \qquad \forall\,i,j = 1,\ldots,n,
\end{equation}
i.e, $P^\top=D_{\boldsymbol{\pi}} P D_{\boldsymbol{\pi}}^{-1}$. With an abuse of notation, the stochastic matrix $P$ is also said to be \textit{reversible}.
\end{definition}
For a reversible Markov chain, we have $D_{\boldsymbol{\pi}}^{-1/2} P^\top D_{\boldsymbol{\pi}}^{1/2} = D_{\boldsymbol{\pi}}^{1/2} P D_{\boldsymbol{\pi}}^{-1/2}$, so $P$ is similar to a symmetric matrix. Consequently, its eigenvalues are real, lie in $[-1,1]$, and $1$ is a simple eigenvalue. 
Moreover, any probability vector $\boldsymbol{\pi}$ satisfying the detailed balance condition~\eqref{eq:detailed_balance} is a stationary distribution. Summing~\eqref{eq:detailed_balance} over $i$ for fixed $j$ gives $\sum_i \pi_i P_{ij} = \sum_i \pi_j P_{ji} = \pi_j \sum_i P_{ji} = \pi_j$; hence, $\boldsymbol{\pi}^\top P = \boldsymbol{\pi}^\top$.
\begin{definition}\label{def:aperiodic-markov-chain}%
An irreducible Markov chain with transition matrix $P$ is said to be \textit{aperiodic} {if there exists an integer $k$ such that, for every pair $i,j$,}  $[P^k]_{ij}>0$. 
\end{definition}
Kemeny's constant is a measure in Markov chain theory that represents  the expected time to reach a random state in a Markov chain, averaged over all possible target states according to the stationary distribution, and it is the same for any starting state. To formally define it, it is necessary to introduce the mean first-passage time $T_{ij}$ from state $ i $ to state $ j $, which is recursively defined~as:
\[
T_{ij} = 
\begin{cases} 
0, & \text{if} \ i = j, \\
1 + \sum_{k \neq i} P_{ik} T_{kj}, & \text{if} \ i \neq j.
\end{cases}
\]
{This quantity represents the expected number of steps to reach $j$ for the first time starting from $i$.}

\begin{definition}\label{def:kemeny_definition}
Let $ P $ be the transition matrix of an irreducible and aperiodic Markov chain, and $ \boldsymbol{\pi}$ be its stationary distribution. \textit{Kemeny's constant} $ \mathcal{K}(P)$ is given by
$
\mathcal{K}(P) = \sum_{j=1}^{n} \pi_j T_{ij}.
$
The value of $\mathcal{K}(P)$ is independent of the starting state~$i$~\cite{doyle2009kemenyconstantmarkovchain}.
\end{definition}
{\begin{remark}
We stress that the assumption of aperiodicity in Definition~\ref{def:kemeny_definition} is not required for the algebraic definition of Kemeny's constant nor for its independence of the initial state. However, we retain it here since, in the irreducible and aperiodic case, the Markov chain is ergodic and admits convergence of $(P^\top)^k$ to $\mathbf{1}\boldsymbol{\pi}^\top$, as $k \rightarrow \infty$, which provides a natural dynamical interpretation of Kemeny's constant in terms of long-run behaviour and mean access times.
\end{remark}}

By using \cite[Theorem 1]{MR636037}, we obtain an expression of Kemeny's constant which allows its computation:
\begin{equation}\label{eq:kemeny-as-trace-with-no-pi}
    \mathcal{K}(P) = \operatorname{tr}\left( (I - P + \mathbf{1}\mathbf{h}^\top )^{-1} \right) - 1, \quad \forall\,\mathbf{h}\in\mathbb{R}^n\,:\,\mathbf{h}^\top\mathbf{1} = 1,
 \end{equation}
{for $I$ the identity matrix, and independently of the vector $\mathbf{h}$.}

In this paper, we will consider stochastic matrices in which  only transitions between selected states are permitted. To handle this setting, we introduce the following notation.

\begin{definition}\label{def:sparsity_pattern}
Given an integer $n\ge 1$, a \emph{pattern} $\mathcal{S}$ is a set of unordered pairs such that
\begin{equation}\label{eq:sparsity_pattern}
    \{i,i\} \in \mathcal{S}\;\forall\,i=1,\ldots,n, \text{ and } \mathcal{S} \subseteq \bigl\{\,\{i,j\} \;\bigm|\; 1 \leq i,j \leq n \,\bigr\}.
\end{equation}
The set \(\mathcal{S}\) describes the admissible positions of nonzero entries in an \(n \times n\) matrix. We say that a matrix \(\Delta \in \mathbb{R}^{n \times n}\) \emph{respects the pattern}~\(\mathcal{S}\) if
\[
\{i,j\} \notin \mathcal{S} \;\;\Rightarrow\;\; \Delta_{ij} = \Delta_{ji} = 0.
\]
The set of all matrices respecting the pattern is denoted by
\[
\mathbb{R}^{n\times n}(\mathcal{S}) \;=\; \bigl\{\, \Delta \in \mathbb{R}^{n\times n} \;\bigm|\; 
\{i,j\} \notin \mathcal{S} \Rightarrow \Delta_{ij} = \Delta_{ji} = 0 \,\bigr\}.
\]
Moreover, we say that a matrix \(\Delta \in \mathbb{R}^{n \times n}\) has the \emph{exact pattern}~\(\mathcal{S}\) if
\[
\Delta_{ij} \neq 0 \text{ and } \Delta_{ji} \neq 0 \;\;\Leftrightarrow\;\; \{i,j\} \in \mathcal{S}.
\]
The set of all such matrices is denoted by
\[
\mathbb{R}^{n\times n}_{\mathrm{exact}}(\mathcal{S}) \;=\; \bigl\{\, \Delta \in \mathbb{R}^{n\times n} \;\bigm|\; 
\Delta_{ij} \neq 0 \text{ and } \Delta_{ji} \neq 0  \;\Leftrightarrow\; \{i,j\} \in \mathcal{S} \,\bigr\}.
\]
Finally, we say that the pattern \(\mathcal{S}\) is \emph{irreducible} if any matrix 
\(\Delta \in \mathbb{R}^{n\times n}_{\mathrm{exact}}(\mathcal{S})\) is irreducible.
\end{definition}

{\begin{remark}
The pattern $\mathcal{S}$ introduced in Definition~\ref{def:sparsity_pattern} admits a natural graph-theoretic interpretation. In particular, it defines an undirected graph on $n$ vertices, where $\{i,j\} \in \mathcal{S}$ corresponds to an edge between nodes $i$ and $j$, and the condition $\{i,i\} \in \mathcal{S}$ for all $i$ corresponds to self-loops at every vertex. In this setting, matrices respecting the pattern have nonzero entries only on the edges of the associated graph. Moreover, the irreducibility of the pattern is equivalent to the connectivity of the underlying undirected graph: a pattern $\mathcal{S}$ is irreducible if and only if the graph induced by $\mathcal{S}$ is connected.
\end{remark}}

\section{Formulating the optimization problems}\label{sec:all-the-problems}

The aim of this section is to formulate the problem of minimizing Kemeny's constant as an optimization problem.

It is worth recalling that, for an irreducible Markov chain with a prescribed stationary probability vector, we know the minimal value of Kemeny’s constant that can be achieved and the form of the matrix that attains it.
\begin{theorem}[{\cite[Theorem~2.2]{MR3218779}}]\label{thm:kirkland_bound}
Suppose that $P$ is an irreducible stochastic matrix of order $n$, and $\boldsymbol{\pi}$ is its stationary distribution, whose components can be assumed to be ordered as $\pi_1 \leq \pi_2 \leq \cdots \leq \pi_n$, {without loss of generality}. Then,
\begin{equation}\label{eq:kirkland-bound}
\mathcal{K}(P) \geq \sum_{j=1}^{n} (j-1) \pi_j.
\end{equation}
Let $P_{n-1}$ be the leading $(n-1) \times (n-1)$ principal submatrix of $P$, and let $\boldsymbol{\pi}_{n-1}$ represent the leading $(n-1)$-subvector of $\boldsymbol{\pi}$. Equality in \eqref{eq:kirkland-bound} holds if and only if  $P_{n-1}$ is nilpotent and
     $\boldsymbol{\pi}_{n-1}^\top (I - P_{n-1})^{-1} \mathbf{1} = \sum_{j=1}^{n} (n-j) \pi_j$.
\end{theorem}
Note that Theorem \ref{thm:kirkland_bound} provides a lower bound for the value of Kemeny's constant, among the set of stochastic matrices with prescribed stationary probability vector $\boldsymbol{\pi}$. %
{As shown in \cite{MR3218779}, when $n\geq 3$, the transition matrices attaining the lower bound in \eqref{eq:kirkland-bound} are far from being reversible. This motivates studying reductions of Kemeny’s constant starting from a prescribed stochastic matrix $P$, rather than seeking global minimizers.
Hence, rather than characterizing the matrices that globally minimize Kemeny’s constant for a prescribed stationary distribution, we study how to modify a given stochastic matrix $P$ as little as possible while achieving a reduction in Kemeny’s constant.}
 This translates into optimizing Kemeny's constant with a (possibly) small perturbation matrix $\Delta$ of the original matrix $P$. Then, while the lower bound in Theorem \ref{thm:kirkland_bound} provides a suitable tool for assessing the reliability of our results, we require that the norm of the perturbation of $P$ should be kept small.
More specifically, we can consider the following optimization problem
\begin{equation}\label{eq:optimization}
\begin{split}
   \min_{\Delta \in \mathbb{R}^{n \times n} } & \operatorname{tr}\left( (I - (P+\Delta) + \mathbf{1}\mathbf{h}^\top)^{-1} 
 \right) + \frac{1}{2}\| \Delta \|_F^2,\\
 \text{ s.t. } &\; \Delta \mathbf{1} = \mathbf{0},
 \;{\boldsymbol{\pi}^\top \Delta = \mathbf{0}},
 \;\text{ and }\;
 \; P + \Delta \geq 0.
\end{split}
\end{equation}
This objective function combines an expression for Kemeny's constant provided in~\eqref{eq:kemeny-as-trace-with-no-pi}, with an additional term equal to half the square of the Frobenius norm of the perturbation $\Delta$. {The latter term is included to penalize large modifications of the original stochastic matrix $P$, since the norm of $\Delta$ can be interpreted as a measure of the effort, energy expenditure, or resources required to implement the perturbation.}
Since the starting matrix $P$ is stochastic, the conditions $\Delta\mathbf{1} = \mathbf{0}$ and $P+\Delta \geq 0$ ensure that $P+\Delta$ also remains stochastic.
{Moreover, the constraint $\boldsymbol{\pi}^\top \Delta = \mathbf{0}$ guarantees that the stationary distribution of the perturbed matrix $P+\Delta$ remains equal to the prescribed stationary distribution $\boldsymbol{\pi}$.}  Unfortunately, the function that maps a stochastic matrix to its Kemeny's constant, as used in the objective of~\eqref{eq:optimization}, is non-convex \cite[\S 3.1]{BulloReview}. 
Note that the objective function in~\eqref{eq:optimization} is a map from $\mathbb{R}^{n\times n}$ to $\mathbb{R} \cup \left\lbrace \infty \right\rbrace $, where the value $\infty$ can be reached by $\Delta$ for which the matrix $P+\Delta$ is reducible. Nevertheless, since the function is lower semicontinuous, level-bounded (for each $\alpha \in \mathbb{R}$, the set of feasible points at which the objective functional does not exceed $\alpha$ is bounded) and proper (there exists at least one feasible point at which it takes a finite value), we apply a generalization of Weierstrass theorem \cite[Theorem 1.9]{rockafellar1998variational}, and obtain that there exists a minimum of the objective function. As a result, even if a {solution to}~\eqref{eq:optimization} exists, it may be generically difficult to numerically estimate it.

Nevertheless, we can construct a convex optimization problem when restricting our analysis to reversible Markov chains (Definition~\ref{def:reversible-markov-chain}), possibly with a fixed pattern, and sharing the stationary distribution $\boldsymbol{\pi}$; {indeed, the restriction on the pattern reflects the fact that, in many applications, only a prescribed set of transition probabilities can be modified, e.g., the ones corresponding to the available or controllable connections in the underlying network.} More specifically, assume that $P$ is a reversible stochastic matrix, that $\boldsymbol{\pi}$ is its stationary distribution and that $P\in\mathbb{R}^{n \times n}(\mathcal{P})$;  Definition~\ref{def:sparsity_pattern}. 
We start considering perturbations $\Delta$ to $P$ such that $\Delta\in\mathbb{R}^{n \times n}(\mathcal{S})$, with $\mathcal{S} \equiv \mathcal{P}$, $P+\Delta$ is stochastic and reversible with pattern $\mathcal{S}\equiv\mathcal{P}$, and  $\boldsymbol{\pi}$ is its stationary distribution. Note that $P+\Delta$ is stochastic, reversible, and such that 
$\boldsymbol{\pi}^\top (P+\Delta)= \boldsymbol{\pi}^\top$, if and only if $\Delta \mathbf{1} = \mathbf{0}$, $P+\Delta\ge 0$, and $ D_{\boldsymbol{\pi}}\Delta =  \Delta^\top D_{\boldsymbol{\pi}}$.

Denote $\boldsymbol{\hat \pi}= \boldsymbol{\pi}^{\nicefrac{1}{2}}$, where the square root is defined component wise. 
By choosing $\mathbf{h}=\boldsymbol{\pi}$ in~\eqref{eq:optimization} and by multiplying the matrix $\left( (I - (P+\Delta) + \mathbf{1}\boldsymbol{\pi}^\top)^{-1} 
 \right)$ to the left by $D_{\boldsymbol{\hat \pi}}$ and to the right by $D_{\boldsymbol{\hat \pi}}^{-1}$, 
 we can construct the following optimization problem:
\begingroup
\allowdisplaybreaks
\begin{equation}\label{eq:optimization-fixed-pi-symmetrized-sparse}
\begin{split}
   \min_{\Delta \in \mathbb{R}^{n\times n}(\mathcal{S})  } & \operatorname{tr}\left( (I - D_{\boldsymbol{\hat \pi}}(P+\Delta)D_{\boldsymbol{\hat \pi}}^{-1} + {\boldsymbol{\hat \pi}}{\boldsymbol{\hat \pi}}^\top )^{-1} 
 \right) + \frac{1}{2}\| \Delta \|_F^2,\\
 \text{ s.t. } &\;  P + \Delta \geq 0,%
 \; D_{\boldsymbol{\pi}}\Delta =  \Delta^\top D_{\boldsymbol{\pi}}, \; \Delta \mathbf{1} = \mathbf{0}. 
\end{split}
\end{equation}

 Then, introducing  the alternative variable $X=P+\Delta$, which is a reversible stochastic matrix with the same pattern $\mathcal{S}$,
the problem \eqref{eq:optimization-fixed-pi-symmetrized-sparse} becomes
\begingroup
\allowdisplaybreaks
\begin{equation}\label{eq:optimization-fixed-pi-symmetrized-in-X}
\begin{split}
   \min_{X \in \mathbb{R}^{n \times n}(\mathcal{S})} & \operatorname{tr}\left( (I - D_{\boldsymbol{\hat\pi}} X D_{\boldsymbol{\hat \pi}}^{-1} + {\boldsymbol{\hat \pi}}{\boldsymbol{\hat \pi}}^\top)^{-1} 
 \right) + \frac{1}{2}\| X-P \|_F^2,\\
 \text{ s.t. } &\; X \geq 0, %
 \;D_{\boldsymbol{\pi}} X = X^{\top} D_{\boldsymbol{\pi}}, 
\; X \mathbf{1} = \mathbf{1}.
\end{split}
\end{equation}
\endgroup

A generalization of this problem can be taken considering the case when the pattern imposed on the perturbation $\Delta$ does not coincide with the pattern $\mathcal{P}$ of the original matrix $P$. 
To this purpose, given a reversible stochastic matrix $P\in\mathbb{R}^{n\times n}(\mathcal{P})$ with stationary distribution $\boldsymbol{\pi}$, and given a pattern $\mathcal{S}$,
we introduce the feasible set 
\begin{equation}
\label{eq:fixed-feasible}
    \mathcal{F}=\{ X \in \mathbb{R}^{n \times n}(\mathcal{S} \cup \mathcal{P}) \,:\, X_{ij}= P_{ij}\mbox{ for } \left\lbrace i,j\right\rbrace \in\mathcal{P}\setminus \mathcal{S}, \, X \mathbf{1} = \mathbf{1}, \, \,D_{\boldsymbol{\pi}}X =  X^\top D_{\boldsymbol{\pi}},\, X \geq 0 \},
\end{equation}
and  the associated optimization problem
\begingroup
\allowdisplaybreaks
\begin{equation}\label{eq:optimization-fixed-pi-fixed}
\begin{split}
   \min_{X \in \mathcal{F}} & \operatorname{tr}\left( (I - D_{\boldsymbol{\hat\pi}} X D_{\boldsymbol{\hat \pi}}^{-1} + {\boldsymbol{\hat \pi}}{\boldsymbol{\hat \pi}}^\top)^{-1} 
 \right) + \frac{1}{2}\| X-P \|_F^2.\\
\end{split}
\end{equation}
\endgroup
Observe that, since $P$ is irreducible, then the pattern $\mathcal P$ is irreducible. Therefore, the pattern $\mathcal S \cup\mathcal P$ is also irreducible. In principle, a matrix $X\in\mathcal F$ might be reducible---note that since the matrix $X$ is reversible, it is equal to a direct sum of two or more stochastic matrices---but in that case, its Kemeny's constant would be infinity, therefore such a matrix cannot minimize the functional in \eqref{eq:optimization-fixed-pi-fixed}.

\subsection{Feasibility and Convexity of the optimization problem}

We now prove that the formulation~\eqref{eq:optimization-fixed-pi-fixed} of the problem is feasible, i.e., the set of constraints $\mathcal{F}$ is not empty. We start proving that, given a prescribed stationary distribution $\boldsymbol{\pi}>0$, there exist reversible stochastic matrices in $\mathbb{R}^{n\times n}(\mathcal{S})$. The following proposition presents conditions on the pattern $\mathcal{S}$ that are necessary to obtain the desired result.

\begin{proposition}\label{pro:pattern-beatrice}
Let $P$ be a reversible stochastic matrix, with stationary distribution $\boldsymbol{\pi}$, such that $P \in \mathbb{R}^{n\times n} (\mathcal{S})$, where $\mathcal{S}$ is defined as in~\eqref{eq:sparsity_pattern}. Then, the set
$
\left\lbrace X \in \mathbb{R}^{n\times n}(\mathcal{S}): X \; \mbox{reversible}, \, \boldsymbol{\pi}^{\top}X= \boldsymbol{\pi}^{\top} \right\rbrace \setminus \left\lbrace P \right\rbrace 
$
is not empty.
\end{proposition}

\begin{proof}
Let $\theta\in\mathbb{R}$ and define $X= P + \theta (P -I)$, where $I$ is the identity matrix. Clearly,  $D_{\boldsymbol{\pi}}X =X^\top D_{\boldsymbol{\pi}}$ and, if $P_{ij}=0$, $i\ne j$, then ${X}_{ij}=0$, therefore $X$ has the pattern of $P+I$. Moreover, by construction, $X\mathbf{1}=\mathbf{1}$. It remains to show that there exist values of $\theta$ such that $X\ge 0$.
If $\left\lbrace i,j\right\rbrace\in\mathcal{S}$ and $i\ne j$, the nonnegativity condition becomes $\theta\ge -1$. If $i=j$, we have $ P_{ii}\ge \theta (1- P_{ii})$, i.e., $\theta\le \frac{P_{ii}}{1-{P}_{ii}}$.
Therefore $X\ge 0$ for any $-1\le\theta\le \min_i {\frac{P_{ii}}{1-P_{ii}}}$.    
\end{proof}
\begin{remark}
    Observe that the same result can be achieved using a convexity argument. Indeed, the matrix $\alpha P + (1-\alpha)I$ belongs to the set for any $\alpha\in (0,1)$.
\end{remark}
More generally, by applying the Metropolis--Hastings adjustment~\cite{Metropolis}, we can prove that the set $\mathcal{F}$ defined in~\eqref{eq:fixed-feasible} is not empty.

\begin{proposition}
\label{pro:metropolis-mod}
    Given a stationary distribution $\boldsymbol{\pi} >0$ and an irreducible pattern $\mathcal{S}$, the set of reversible stochastic matrices $X\in \mathbb{R}^{n\times n}(\mathcal{S})$, with pattern $\mathcal{S}$ as in~\eqref{eq:sparsity_pattern} is not empty.

    Moreover, given a reversible stochastic  matrix $P \in \mathbb{R}^{n \times n}(\mathcal{P})$ with pattern $\mathcal{P}$ and  stationary distribution $\boldsymbol{\pi}$, and given a pattern $\mathcal{S}$, then the set $\mathcal{F}$ in~\eqref{eq:fixed-feasible} is not empty.
\end{proposition}

\begin{proof}
For the first part of the proof, consider the matrix $A$ associated with the pattern $\mathcal{S}$ by
\begin{equation}\label{eq:adjacency}
    A_{ij} = \begin{cases}
    1, & \{i,j\} \in \mathcal{S},\\
    0, & \text{otherwise},
\end{cases}
\end{equation}
and define the stochastic matrix $Q = D_{A\mathbf{1}}^{-1} A$, which, by construction, has the pattern $\mathcal{S}$.

Define the transition matrix $X$ with the same pattern $\mathcal{S}$ by applying the Metropolis--Hastings adjustment to~$Q$: We set, for each $i \neq j$
\[
X_{ij} = \begin{cases} \displaystyle Q_{ij} \cdot \min\left\lbrace 1, \frac{\pi_j Q_{ji}}{\pi_i Q_{ij}} \right\rbrace, & Q_{ij} \neq 0,\\[1em]
0,& \text{otherwise},
\end{cases}
\]
and define the diagonal by $X_{ii} = 1 - \sum_{j \neq i} X_{ij}$ to ensure row-stochasticity. %
The detailed balance condition~\eqref{eq:detailed_balance} holds by the pattern symmetry of the matrix $Q$, and the construction preserves the pattern $\mathcal{S}$ by only modifying entries where $Q_{ij} > 0$, and possibly the diagonal entries to enforce the condition on the sum over the rows. {Since the pattern $\mathcal{S}$ is irreducible, the associated undirected graph is connected. By construction, the matrix $X$ has strictly positive off-diagonal entries for all the indices $\{i,j\} \in \mathcal{S}$ (possibly after Metropolis adjustment and scaling), and no other off-diagonal positive entries. Hence, the directed graph induced by the support of $X$ coincides with the connected graph associated with $\mathcal{S}$; therefore, $X$ is irreducible as well.}

For the second part of the proposition, we construct a reversible Markov chain whose transition matrix belongs to $\mathcal{F}$ as follows: given a stochastic $Q \in \mathbb{R}^{n\times n}_{\rm{exact}}(\mathcal{S})$---possibly built as $Q = D_{A\mathbf{1}}^{-1}A$ with $A$ as in~\eqref{eq:adjacency}---we set for~$i\neq j$
\begin{align*}
    X_{ij} = \left\lbrace \begin{array}{l}
     P_{ij}, \quad \mbox{for } \left\lbrace i,j\right\rbrace\not\in \mathcal{S}, \\
    \delta\cdot Q_{ij}\cdot \min\left\lbrace 1,\frac{\pi_j Q_{ji}}{\pi_i Q_{ij}}\right\rbrace, \quad \mbox{for } \left\lbrace i,j\right\rbrace \in \mathcal{S},
    \end{array}\right.
\end{align*}
where we choose 
{$\delta \in \left[ 0,1 -\|P^0\mathbf{1}\|_{\infty} \right]$}, for
\begin{align}
\label{eq:P0}
    P^0_{ij} = \left\lbrace \begin{array}{l}
     P_{ij}, \quad \mbox{for } \left\lbrace i,j\right\rbrace\not\in \mathcal{S}, \\
    0, \quad \mbox{for } \left\lbrace i,j\right\rbrace \in \mathcal{S}.
    \end{array}\right.
\end{align}
We define the diagonal entries as $X_{ii}=1-\sum_{j \neq i} X_{ij} $, to ensure row-stochasticity. The balance condition may be proved entrywise, for entries $\left\lbrace i,j\right\rbrace \in \mathcal{S}$, while it holds automatically for entries $\left\lbrace i,j\right\rbrace\not\in \mathcal{S}$, {and irreducibility of $X$ follows from irreducibility of the pattern $\mathcal{S}$}.
\end{proof}

\begin{remark}
\label{remark: only one}
    Observe that if the pattern $\mathcal{S}$, in the hypotheses of the second part of the Proposition~\ref{pro:metropolis-mod}, is such that $\|P^0\mathbf{1}\|_{\infty} =1$, then the set $\mathcal{F}$ in~\eqref{eq:fixed-feasible} consists only of the matrix $P$, since $\delta = 0$. In such case, the optimization problem~\eqref{eq:optimization-fixed-pi-fixed} becomes meaningless. Moreover, the construction in Propositions~\ref{pro:pattern-beatrice} and~\ref{pro:metropolis-mod} holds true for matrices $X \in \mathbb{R}^{n\times n}_{\rm exact}(\mathcal{S})$.
\end{remark}

The need for including the $\{i,i\}$ pairs in the pattern $\mathcal{S}$ in Definition~\ref{def:sparsity_pattern} is now clear. 
Indeed, firstly we construct a stochastic matrix $Q$, with the desired pattern $\mathcal{S}$, but with a stationary distribution different from $\boldsymbol{\pi}$. Secondly, we adjust the entries obtaining a reversible matrix with both pattern $\mathcal{S}$ and stationary distribution $\boldsymbol{\pi}$. To this end, the presence of the diagonal entries and the possibility to adjust them was essential for the construction of the matrix $X$.
This condition is the same as that used in~\cite{MR4822696} for assigning a target stationary distribution to a given Markov chain.

We now prove that the formulation~\eqref{eq:optimization-fixed-pi-fixed} has a unique solution through a convexity argument that follows the proof of the analogous result~\cite[Theorem~2]{Bullo}. The convexity results for the problems~\eqref{eq:optimization-fixed-pi-symmetrized-sparse} and \eqref{eq:optimization-fixed-pi-symmetrized-in-X} will then directly follow, by choosing $\mathcal{S} \equiv \mathcal{P}$.

\begin{proposition}\label{pro:the_problem_is_convex}
Given a reversible stochastic matrix $P \in \mathbb{R}^{n \times n}(\mathcal{P})$, with stationary distribution $\boldsymbol{\pi}$, and given a pattern $\mathcal{S}$, then the set $\mathcal{F}$ in~\eqref{eq:fixed-feasible} is convex. Furthermore, the function $f$ 
    \[
    \begin{split}
            X \mapsto f(X) = &\; \operatorname{tr}\left( (I - D_{\boldsymbol{\hat \pi}} X D_{\boldsymbol{\hat \pi}}^{-1} + {\boldsymbol{\hat \pi}}{\boldsymbol{\hat \pi}}^{\top})^{-1}  
 \right) %
 + \frac{1}{2}\| X-P \|_F^2,
    \end{split}
    \]
    is convex over $\mathcal{F}$.
\end{proposition}

\begin{proof}
Let $X_1,X_2 \in \mathcal{F}$, then for all $\theta \in [0,1]$, by direct inspection, we have
\[
\begin{split}
 (\theta X_1 + (1-\theta) X_2) \mathbf{1} = &\; \theta X_1 \mathbf{1} + (1-\theta) X_2 \mathbf{1} = \mathbf{1},\\
 D_{\boldsymbol{\pi}} (\theta X_1 + (1-\theta) X_2) = & \; (\theta X_1 + (1-\theta) X_2)^\top D_{\boldsymbol{\pi}}, \\
 P + (\theta X_1 + (1-\theta) X_2) = &\; \theta (P + X_1) + (1-\theta) (P+ X_2) \geq 0,
\end{split}
\]
hence the set is convex, since $(\theta X_1+ (1-\theta)X_2) \in \mathbb{R}^{n\times n}(\mathcal{S}\cup \mathcal{P})$ and $\theta{\left(X_1\right)}_{ij} + (1-\theta) {\left(X_2\right)}_{ij} = P_{ij}$ for each entry in $\mathcal{P}\setminus \mathcal{S}$. For the convexity of the function, we observe that $W = I - D_{\boldsymbol{\hat \pi}}X D_{\boldsymbol{\hat \pi}}^{-1} + {\boldsymbol{\hat \pi}} {\boldsymbol{\hat \pi}}^{\top}$ is a symmetric positive definite matrix by construction.

Indeed, since the matrix $D_{\boldsymbol{\hat \pi}}X D_{\boldsymbol{\hat \pi}}^{-1}$ is similar to $X$, we get that its spectral radius is equal to $1$. Moreover, since $X\in \mathcal{F}$, the matrix $D_{\boldsymbol{\hat \pi}}X D_{\boldsymbol{\hat \pi}}^{-1}$ is symmetric. We get that the eigenvalues of $I - D_{\boldsymbol{\hat \pi}}X D_{\boldsymbol{\hat \pi}}^{-1}$ are real and nonnegative. Since $\boldsymbol{\hat \pi}$ is in the kernel of the matrix $I - D_{\boldsymbol{\hat \pi}}X D_{\boldsymbol{\hat \pi}}^{-1}$ and 
$
W\boldsymbol{\hat \pi} = \boldsymbol{\hat \pi} - D_{\boldsymbol{\hat \pi}} X \mathbf{1}  + (\boldsymbol{\hat \pi}^{\top}\boldsymbol{\hat \pi}) \boldsymbol{\hat \pi} = (\boldsymbol{\hat \pi}^{\top}\boldsymbol{\hat \pi}) \boldsymbol{\hat \pi}$,
we conclude that $W$ is symmetric positive definite.

The function %
$\operatorname{tr}\left( (I - D_{\boldsymbol{\hat \pi}} X D_{\boldsymbol{\hat \pi}}^{-1} + {\boldsymbol{\hat \pi}}{\boldsymbol{\hat \pi}}^{\top})^{-1}\right)$ is then the composition of the function $Y\mapsto\operatorname{tr}(Y^{-1})$ and an affine mapping from the convex set $\mathcal{F}$ to a subset of symmetric positive definite matrix---the inverse preserve the sign of the eigenvalues---hence it is convex~\cite[\S 3.2.2]{ConvexOptimization}. Then $f(X)$ is convex as the sum of two convex functions.
\end{proof}

\section{Formulation as constrained optimization}\label{sec:solving-all-the-problems}

We solve the optimization problems~\eqref{eq:optimization-fixed-pi-symmetrized-sparse} by formulating it as a constrained minimization problem in the following general form:
\begin{equation}\label{eq:optimization_general_form}
\begin{split}
\min_{\boldsymbol{\delta} \in \mathbb{R}^{m}} \; & g(\boldsymbol{\delta})\,:\,\mathbb{R}^{m} \rightarrow \mathbb{R},\\
\text{s.t. } &\; \begin{array}{ll}
A \boldsymbol{\delta} = \mathbf{b}, & \text{(equality constraints)},\\
\mathbf{l}_b \leq \boldsymbol{\delta}, & \text{(bound constraints)}.
\end{array}
\end{split}    
\end{equation}
Here, \( g(\boldsymbol{\delta}) \) represents the objective function to be minimized, \( A \boldsymbol{\delta} = \mathbf{b} \) denotes the set of linear equality constraints, and the additional lower bound constraints on \( \boldsymbol{\delta} \) are represented by \( \mathbf{l}_b \leq \boldsymbol{\delta} \).

To express $A$, $\mathbf{b}$, and $\mathbf{l}_b$ in~\eqref{eq:optimization_general_form} in our case, we introduce two operations: $\operatorname{vec}(\cdot)$ and $\operatorname{mat}(\cdot)$. The operator $\operatorname{vec}(\cdot)$ takes a matrix of size \( n \times n \) and stacks its columns into a vector of size \( n^2 \), while $\operatorname{mat}(\cdot)$ reverses this operation, mapping a vector of size \( n^2 \) back to a matrix of size \( n \times n \). These operations satisfy the property \( \operatorname{mat}(\operatorname{vec}(A)) = A \) for any matrix \( A \in \mathbb{R}^{n \times n} \).

We can now reformulate the optimization problems on the matrix \( \Delta \) as problems involving the vector variable \( \boldsymbol{\delta} = \operatorname{vec}(\Delta) \). As a constructive tool, we write the vector formulation for~\eqref{eq:optimization-fixed-pi-symmetrized-sparse}, choosing $\mathcal{S}= \bigl\{\,\{i,j\} \;\bigm|\; 1 \leq i,j \leq n \,\bigr\}$:
\begin{equation}\label{eq:optimization-fixed-pi-symmetrized-vectorized}
\begin{split}
   \min_{\boldsymbol{\delta} \in  \mathbb{R}^{n^2}} & \; g(\boldsymbol{\delta}) \equiv \operatorname{tr}\left( (I - D_{\boldsymbol{\hat \pi}} (P+\operatorname{mat}(\boldsymbol{\delta}))D_{\boldsymbol{\hat \pi}}^{-1} + \boldsymbol{\hat \pi}\boldsymbol{\hat \pi}^\top)^{-1} 
 \right) + \frac{1}{2}\| \boldsymbol{\delta} \|_F^2%
 \,:\,\mathbb{R}^{n^2} \rightarrow \mathbb{R},\\
 \text{s.t. } & \; A \boldsymbol{\delta} \equiv \begin{bmatrix}
  \mathbf{1}^\top \otimes I  \\
  I \otimes D_{\boldsymbol{\pi}} - (D_{\boldsymbol{\pi}} \otimes I ) T
 \end{bmatrix}
 \boldsymbol{\delta} = \mathbf{0} \equiv \mathbf{b},%
 \; \mathbf{l}_b \equiv -\operatorname{vec}(P) \leq \boldsymbol{\delta},
\end{split}
\end{equation}
where \( I \) is the identity matrix, $\otimes$ denotes the Kronecker product of two matrices, and $T$ is the permutation matrix, oftentimes called the orthogonal commutation matrix, mapping the stacking $\operatorname{vec}(\Delta^\top)$ into $\boldsymbol{\delta}$. %

IPMs are effective for problems of the form~\eqref{eq:optimization_general_form}~\cite{MR4865731}. They handle inequality constraints via barrier functions, keeping iterates strictly feasible, and solve a sequence of approximations to converge to the optimum. IPMs are particularly efficient for large-scale problems, offering polynomial-time complexity and high accuracy. To apply them to~\eqref{eq:optimization-fixed-pi-symmetrized-vectorized}, we require the gradient of $g(\boldsymbol{\delta})$, and generally its Hessian. It is convenient to express these derivatives in matrix form, as summarized in the following technical lemma.
\begin{lemma}\label{lem:boring_computations}
Given the function
$
{\psi}(X) = \operatorname{tr}( (\Omega - D_{\mathbf{v}} X D_{\mathbf{w}})^{-1} ) + \frac{1}{2} \| X \|_F^2 \,:\,\mathbb{R}^{n \times n} \rightarrow \mathbb{R},
$
for $\Omega,D_{\mathbf{v}},D_{\mathbf{w}} \in \mathbb{R}^{n \times n}$ and 
$\mathbf{v},\mathbf{w}\in \mathbb{R}^n$, the entries of the gradient of ${\psi}(X)$ are given by
\[
\frac{\partial {\psi}(X)}{\partial X_{ij}} = v_i w_j \, \mathbf{e}_j^\top ( \Omega - D_{\mathbf{v}} X D_{\mathbf{w}} )^{-2} \mathbf{e}_i + X_{ij},
\]
with $\mathbf{e}_j$, $\mathbf{e}_i$ the $j$th and $i$th vector of the canonical basis of $\mathbb{R}^n$.
The entries of the Hessian of ${\psi}(X)$ are expressed as
\[
\begin{split}
\frac{\partial^2 {\psi}(X)}{\partial X_{ij}\partial X_{hk}} = & v_i w_j v_h w_k \Big( \mathbf{e}_j^\top (\Omega - D_{\mathbf{v}} X D_{\mathbf{w}}) ^{-1}\mathbf{e}_h \; \mathbf{e}_k^\top (\Omega - D_{\mathbf{v}} X D_{\mathbf{w}}) ^{-2}\mathbf{e}_i + \\ & +  
\mathbf{e}_j^\top (\Omega - D_{\mathbf{v}} X D_{\mathbf{w}}) ^{-2}\mathbf{e}_h \; \mathbf{e}_k^\top (\Omega - D_{\mathbf{v}} X D_{\mathbf{w}}) ^{-1}\mathbf{e}_i \Big) + \delta_{ih}\delta_{jk},
\end{split}
\]
with $\delta_{rs}$ the Kronecker delta, with $r,s=1,\ldots,n$.
\end{lemma}

\begin{proof}
    The formulae for $\frac{\partial {\psi}(X)}{\partial X_{ij}}$ and $\frac{\partial^2 {\psi}(X)}{\partial X_{ij} \;\partial X_{hk}}$ {are derived} from direct computation, exploiting the linearity and the cyclic property of the trace.
\end{proof}

Then, utilizing the expression of $g(\boldsymbol{\delta})$ in~\eqref{eq:optimization-fixed-pi-symmetrized-vectorized}, we have the following formulae for the gradient and the Hessian.

\begin{proposition}\label{pro:simple-gradient}
The gradient for $g(\boldsymbol{\delta})$ in~\eqref{eq:optimization-fixed-pi-symmetrized-vectorized} is given by
\[
\nabla\,g(\boldsymbol{\delta}) = \operatorname{vec}\left( \hat{\Pi} \odot \left( I - D_{\boldsymbol{\hat \pi}} (P+\operatorname{mat}(\boldsymbol{\delta}))D_{\boldsymbol{\hat \pi}}^{-1} + {\boldsymbol{\hat \pi}}{\boldsymbol{\hat \pi}}^\top \right)^{-2} \right) + \boldsymbol{\delta},\]
for $\hat{\Pi}$ the scaling matrix with entries $\hat{\Pi}_{ij} = \nicefrac{\hat{\pi}_i}{\hat{\pi}_j}$, and $\odot$ the entrywise product.

The entries of the Hessian of the function ${g}(\boldsymbol{\delta})$ are given by 
\[
\begin{split}
    \mathbf{H}_{p,q} = & \frac{\hat{\pi}_i}{\hat{\pi}_j} \frac{\hat{\pi}_h}{\hat{\pi}_k} \left( \mathbf{e}_j^\top H_s(\Delta) ^{-1}\mathbf{e}_h\; \mathbf{e}_k^\top H_s(\Delta) ^{-2}\mathbf{e}_i  +
\mathbf{e}_j^\top H_s(\Delta) ^{-2}\mathbf{e}_h\; \mathbf{e}_k^\top H_s(\Delta) ^{-1}\mathbf{e}_i  \right) + \delta_{pq},
\end{split}
\]
for $p = i + (j-1)n$, $q = h + (k-1)n$ and $i,j,h,k=1,\ldots,n$, $\delta_{pq}$ the Kronecker delta, and $H_s(\Delta) = H_s(\operatorname{mat}(\boldsymbol{\delta}))$ defined by
$
H_s(\Delta)=I - D_{\boldsymbol{\hat\pi}}(P+\Delta)D_{\boldsymbol{\hat\pi}}^{-1} + {\boldsymbol{\hat\pi}}{\boldsymbol{\hat\pi}}^\top.
$
\end{proposition}

\begin{proof}
The formulae follow by Lemma~\ref{lem:boring_computations} setting $\Omega = I - D_{\boldsymbol{\hat \pi}} P D_{\boldsymbol{\hat \pi}}^{-1} + {\boldsymbol{\hat\pi}}{\boldsymbol{\hat\pi}}^\top$, $D_{\mathbf{v}} = D_{\boldsymbol{\hat \pi}}$ and $D_{\mathbf{w}} = D_{\boldsymbol{\hat \pi}}^{-1}$.
\end{proof}

The problem \eqref{eq:optimization-fixed-pi-symmetrized-sparse} can be solved directly in MATLAB using \texttt{fmincon}. Our \texttt{GitHub} implementation provides exact expressions for the objective, gradient, and Hessian (see Proposition~\ref{pro:simple-gradient}). Since forming the Hessian may be costly, we accelerate the computation by implementing core routines as \texttt{mex} files in~C.

Since we are solving a constrained optimization problem, the optimality conditions for the solutions can be investigated via the KKT conditions. Having proved that the minimization problem in~\eqref{eq:optimization-fixed-pi-fixed} is convex in Proposition~\ref{pro:the_problem_is_convex}, we get that the KKT conditions are sufficient and necessary conditions for global optimality~\cite[Section 5.5.3]{ConvexOptimization}. Then, the minimizers of problem~\eqref{eq:optimization-fixed-pi-fixed} coincide with its KKT points. The solutions of the constrained problem are computed via MATLAB \texttt{fmincon}. It employs an IPM solver via a primal-dual log-barrier method~\cite{Byrd}, which iteratively solves a sequence of perturbed KKT systems and is guaranteed to converge to a KKT point of~\eqref{eq:optimization-fixed-pi-fixed}.

\subsection{Adding the constraints on the pattern of the perturbation}
\label{subsec:euclidean_optimizer_sub pattern}
Let us now consider the case where the pattern $\mathcal{S}$ chosen on the perturbation matrix $\Delta$ is not necessary equal to the pattern $\mathcal{P}$ of the starting matrix $P$. In this case, we can represent the matrix \( \Delta  \in \mathbb{R}^{n\times n}(\mathcal{S})\) using a coordinate format, denoted as the triple of vectors \( (\mathbf{r}, \mathbf{c}, \boldsymbol{\delta}) \). These three vectors respectively contain the row indices, column indices, and the possible nonzero entries and have size $m= \left| \mathcal{S} \right|$. We then define the function \( \operatorname{mat}(\mathbf{r}, \mathbf{c}, \boldsymbol{\delta}) \), which constructs the matrix from this coordinate representation. Additionally, the function \( \operatorname{vec}(\mathbf{r},\mathbf{c},\Delta) \) returns only the entries indexed by $\mathbf{r}$ and $\mathbf{c}$.
Formally, given the pattern $\mathcal{S} = \{(r_k,c_k)\}_{k=1}^m$, we construct the vectors $\mathbf{r}= (r_1,\ldots, r_m)$, $\mathbf{c}= (c_1,\ldots, c_m)$. For any matrix $\Delta \in \mathbb{R}^{n\times n}(\mathcal{S})$, we define $
\operatorname{vec}(\mathbf r,\mathbf c,\Delta)
=
\boldsymbol{\delta}\in\mathbb R^m$, such that $\delta_k=\Delta_{r_k,c_k}$, for $k=1,\dots,m.$ Conversely, given a vector $\boldsymbol{\delta} \in \mathbb{R}^m$, we may construct the matrix $\operatorname{mat}(\mathbf r,\mathbf c,\boldsymbol{\delta})
=
\Delta\in\mathbb R^{n\times n}(\mathcal S)$, such that $\Delta_{r_k,c_k}=\delta_k$, for $k=1,\ldots,m$, and zero elsewhere.

We may then write the vector formulation in \eqref{eq:optimization-fixed-pi-symmetrized-vectorized}, with this modification, as:
\begin{equation}
\label{eq:constrained-optimizer-sparse}
g(\boldsymbol{\delta}) = \operatorname{tr}\left( \left(I - D_{\boldsymbol{\hat \pi}}\left(P + \operatorname{mat}(\mathbf{r}, \mathbf{c}, \boldsymbol{\delta})\right)D_{\boldsymbol{\hat \pi}}^{-1} + \boldsymbol{\hat \pi}\boldsymbol{\hat \pi}^\top\right)^{-1} \right) + \frac{1}{2} \|\boldsymbol{\delta}\|_F^2\,:\,\mathbb{R}^m \rightarrow \mathbb{R}.
\end{equation}
To complete the IPM formulation in~\eqref{eq:optimization-fixed-pi-symmetrized-vectorized}, we construct a vector of inequality constraints \( \mathbf{l}_b \) equal to \( \mathbf{l}_b = - \operatorname{vec}(\mathbf{r}, \mathbf{c}, P) \), while to formulate the equality constraints, we define the matrix \( \Gamma \in \mathbb{R}^{m \times n^2} \), such that \( \Gamma \mathbf{v} = \operatorname{vec}(\mathbf{r}, \mathbf{c}, \operatorname{mat}(\mathbf{v})) \). Observe that the matrix $\Gamma$ acts as a projector, restricting a generic vector in \( \mathbb{R}^{n^2} \) to the entries corresponding to the pattern \( \mathcal{S} \). Consequently, in the set of equality constraints, the matrix \( A \) in~\eqref{eq:optimization-fixed-pi-symmetrized-vectorized} is replaced by \( A \begin{bmatrix}\Gamma & \Gamma \end{bmatrix}^\top \), while the vector \( \mathbf{b} \) becomes the zero vector of size \( n^2 + n\). %

The expression of the gradient is then a direct consequence of the result in Proposition~\ref{pro:simple-gradient}.
\begin{corollary}
The gradient of the objective function in~\eqref{eq:constrained-optimizer-sparse}, with $\Delta = \operatorname{mat}(\mathbf{r},\mathbf{c},\boldsymbol{\delta})$ a matrix with pattern $\mathcal{S}$ is given~by
\[
\nabla\,g(\boldsymbol{\delta}) = \operatorname{vec}\left(\mathbf{r},\mathbf{c}, \hat{\Pi} \odot \left( I - D_{\boldsymbol{\hat \pi}}(P+\operatorname{mat}(\mathbf{r},\mathbf{c},\boldsymbol{\delta}))D_{\boldsymbol{\hat \pi}}^{-1} + {\boldsymbol{\hat \pi}}{\boldsymbol{\hat \pi}}^\top \right)^{-2} \right) + \boldsymbol{\delta},\]
for $\hat{\Pi}$ the scaling matrix with entries $\hat{\Pi}_{ij} = \nicefrac{\hat{\pi}_i}{\hat{\pi}_j}$, and $\odot$ the entrywise product.
\end{corollary}
The expression of the Hessian can then be obtained again as a corollary to Proposition~\ref{pro:simple-gradient} by restricting the set of indices $i,j,h,k$ from $1,\ldots,n$ only the values corresponding to the entries in the pattern~$\mathcal{S}$.

%
%
\section{A Riemannian optimization based approach}\label{sec:riemannian_formulation}
As we have anticipated, beyond the constrained approach of Section~\ref{sec:solving-all-the-problems}, problem~\eqref{eq:optimization-fixed-pi-symmetrized-in-X} can also be addressed via Riemannian optimization. {We refer the reader to \cite{BoumalBook, AbsilBook} for introductory books on the topic.} We therefore examine the feasible set and its manifold structure. Assume the initial stochastic matrix $P$ is reversible with stationary distribution $\boldsymbol{\pi}$ and pattern $\mathcal{P}$ as in~\eqref{eq:sparsity_pattern}. We first set $\mathcal{S}\equiv\mathcal{P}$ and consider variable matrices $X \in \mathbb{R}^{n\times n}_{\mathrm{exact}}(\mathcal{S})$ with strictly positive entries on $\{i,j\}\in\mathcal{S}$. Note that $P\in\mathbb{R}^{n\times n}(\mathcal{P})$ may contain zeros on $\mathcal{P}$, unlike~$X$.
More precisely, we are choosing the non-empty feasible set (Remark~\ref{remark: only one}):
\begin{equation*}
    \mathcal{M} = \left\lbrace X \in \mathbb{R}^{n\times n}_{\text{exact}}(\mathcal{S}) : X \; \mbox{stochastic and reversible}, D_{\boldsymbol{\pi}}X = X^{\top} D_{\boldsymbol{\pi}} \right\rbrace
\end{equation*}
and set as a variable $X= P+\Delta$. %
We may then rewrite~\eqref{eq:optimization-fixed-pi-symmetrized-in-X} as
\begin{align*}
    \min_{X \in \mathcal{M}} \operatorname{tr}\left( (I - D_{\boldsymbol{\hat\pi}} X D_{\boldsymbol{\hat \pi}}^{-1} + {\boldsymbol{\hat \pi}}{\boldsymbol{\hat \pi}}^\top)^{-1} 
 \right) + \frac{1}{2}\| X-P \|_F^2,
\end{align*}
where we recall that $\hat{\boldsymbol{\pi}} = \boldsymbol{\pi}^{1/2}$.
Following the idea proposed in~\cite{Reversible}, we introduce the set
\begin{equation}
\label{eq:manifold_sparse_reversible}
\mathcal{M}_{\boldsymbol{\pi}} = \left\lbrace X \in \mathbb{R}^{n \times n}_{\text{exact}}(\mathcal{S}) : X = X^\top, X \hat{\boldsymbol{\pi}} = \hat{\boldsymbol{\pi}}, X_{ij} > 0 \text{ for } \left\lbrace i,j\right\rbrace \in \mathcal{S} \right\rbrace,
\end{equation}
whose construction is justified by a change of variable and a manipulation of the balance condition \eqref{eq:detailed_balance}:
\begin{gather*}
  D_{\boldsymbol{\pi}} (P + \Delta) = (P + \Delta)^{\top} D_{\boldsymbol{\pi}} \iff   D_{\boldsymbol{\pi}} (P+\Delta) D_{\boldsymbol{\pi}}^{-1} = (P+\Delta)^{\top} \\ \iff 
  \overbrace{D_{\boldsymbol{\hat \pi}} \underbrace{(P+ \Delta)}_{X \in \mathcal{M}} D_{\boldsymbol{\hat \pi}}^{-1}}^{X \in \mathcal{M}_{\boldsymbol{\pi}}} = \overbrace{D_{\boldsymbol{\hat \pi}}^{-1} \underbrace{(P+\Delta)^{\top}}_{X^{\top} \in \mathcal{M}} D_{\boldsymbol{\hat \pi}}}^{X^\top \in \mathcal{M}_{\boldsymbol{\pi}}},
\end{gather*}
and the stochasticity condition leads to $X\boldsymbol{\hat\pi} = D_{\boldsymbol{\hat \pi}} (P + \Delta) \mathbf{1} = \boldsymbol{\hat \pi}$, with $X \in \mathcal{M}_{\boldsymbol{\pi}}$.
By this change, the original minimization problem on $\mathcal{M}$ can be reformulated as $\min_{X \in \mathcal{M}_{\boldsymbol{\pi}}} f(X),
$ that is an optimization problem over the set~\eqref{eq:manifold_sparse_reversible}, 
where the functional $f(X)$ is defined as 
   \begin{equation}
   \label{eq:functional_min_X}
   f(X)=  \operatorname{tr} \left( (I - X + {\boldsymbol{\hat \pi}}{\boldsymbol{\hat \pi}}^{\top})^{-1} \right) + \frac{1}{2} \| D_{\boldsymbol{\hat \pi}}^{-1} X D_{\boldsymbol{\hat \pi}}  - P\|_F^2.
   \end{equation}

The set $\mathcal{M}_{\boldsymbol{\pi}}$ is a manifold of symmetrized reversible stochastic matrices satisfying the linear constraints imposed by $\mathcal{S}$; see~\cite[\S3]{Reversible}. By Proposition~\ref{pro:pattern-beatrice} or the Metropolis–Hastings construction in Proposition~\ref{pro:metropolis-mod} (see also Remark~\ref{remark: only one}), this manifold is non-empty, ensuring the existence of feasible points and enabling a manifold-based optimization strategy~\cite{AbsilBook}.
To this end, we derive the differential and Riemannian structure of $\mathcal{M}_{\boldsymbol{\pi}}$ in~\eqref{eq:manifold_sparse_reversible}. We rely on standard notions (tangent space, retraction) as in~\cite[\S3,\S5]{BoumalBook}, and restrict attention to first-order geometry—eschewing the Riemannian Hessian and Levi–Civita connection. The resulting framework supports gradient-based methods, which we will see are sufficient for our purposes.

Unlike the IPM previously described, Riemannian optimizers work directly on the manifold, creating a sequence of iterates that belong to the feasible set $\mathcal{M}_{\boldsymbol{\pi}}$ by construction. In this setting, the appropriate notion of optimality is not given by the KKT conditions, but rather by the first-order optimality condition on the manifold, namely the vanishing of the Riemannian gradient. In particular, the Riemannian conjugate gradient method is guaranteed to produce a sequence whose accumulation points satisfy this condition~\cite[Theorem 4.3.1]{AbsilBook}.

The structure of the tangent space is inherited from \cite[Lemma 3.2]{Reversible}, by adding the constraint on the pattern~$\mathcal{S}$. 
\begin{corollary}
   The tangent space at a point $X \in \mathcal{M}_{\boldsymbol{\pi}}$ is
\begin{equation}\label{eq:tangent_space}
   \mathcal{T}_{X} \mathcal{M}_{\boldsymbol{\pi}} = \left\lbrace \xi_X \in \mathbb{R}^{n \times n}(\mathcal{S}) : \xi_X = \xi_X^\top, \xi_X \hat{\boldsymbol{\pi}} =0 %
   \right\rbrace.
\end{equation}
\end{corollary}
To obtain a Riemannian manifold we need to endow the tangent space~\eqref{eq:tangent_space} with an inner product. Note that in this setting, the choice of the Fisher metric made in~\cite{Reversible}, or in~\cite{Douik} for the \emph{multinomial manifold}, does not define an inner product on the tangent space, due to the presence of zero entries in the matrices. To circumvent this limitation, we consider a variant of this metric, defined as follows: for $\xi_X,\eta_X \in \mathcal{T}_X \mathcal{M}_{\boldsymbol{\pi}}$
\begin{equation}
\label{eq:fisher_sparse}
   \left\langle \xi_X,\eta_X \right\rangle_X = \sum_{X_{ij} \neq 0} \frac{(\xi_X)_{ij} (\eta_X)_{ij}}{X_{ij}}  =  
   \operatorname{tr}\left( ( \xi_X \odiv X ) \eta_X^\top \right),
\end{equation}
where we denote by $\odiv$ the following operation: given $A \in \mathbb{R}^{n \times n}(\mathcal{S}), B \in \mathbb{R}^{n \times n}_{\text{exact}}(\mathcal{S})$, we define
\begin{align*}
   (A \odiv B)_{ij} = \left\lbrace \begin{array}{c}
       \frac{A_{ij}}{B_{ij}}, \quad \mbox{for } B_{ij} \neq 0,  \\
       0, \quad \; \mbox{for } B_{ij} = 0.
   \end{array} \right.
\end{align*}
This choice represents an inner product on the tangent space; indeed, by direct computation it can be shown that~\eqref{eq:fisher_sparse} is bi-linear, symmetric and positive definite on $\mathcal{T}_X\mathcal{M}_{\boldsymbol{\pi}}$.

\begin{lemma}
   The function \eqref{eq:fisher_sparse} is an inner product on the tangent space $\mathcal{T}_X\mathcal{M}_{\boldsymbol{\pi}}$ in~\eqref{eq:tangent_space}.
\end{lemma}

%
%
%

Hence, this choice of inner product \eqref{eq:fisher_sparse} defines a metric on the tangent space. We only need to check that this is a Riemannian metric~\cite[Deﬁnition 3.52]{BoumalBook}, that is for all smooth vectors fields $V,W$ on $\mathcal{M}_{\boldsymbol{\pi}}$ the function $X \mapsto \left\langle V(X), W(X) \right\rangle_X$ is \emph{smooth}---where a smooth vector field on the manifold $\mathcal{M}_{\boldsymbol{\pi}}$ is a smooth map $V: \mathcal{M}_{\boldsymbol{\pi}} \mapsto \mathcal{T}\mathcal{M}_{\boldsymbol{\pi}}$ such that $V(X) \in \mathcal{T}_X \mathcal{M}_{\boldsymbol{\pi}}$, for every $X \in \mathcal{M}_{\boldsymbol{\pi}}$. This holds observing that the inner product is a composition of smooth functions, since $X \in \mathbb{R}^{n\times n}_{\rm exact}(\mathcal{S})$. 

\begin{remark}
   It is worth noticing that the manifold $\mathcal{M}_{\boldsymbol{\pi}}$ in \eqref{eq:manifold_sparse_reversible} is not a submanifold of the manifold 
   $\left\lbrace X \in \mathbb{R}^{n\times  n}: X= X^{\top}, X \boldsymbol{\hat \pi} = \boldsymbol{\hat \pi}, \, X >  0 \right\rbrace$,
   studied in \cite{Reversible}. Indeed, the elements in the manifold $\mathcal{M}_{\boldsymbol{\pi}}$ do not satisfy the entry-wise positivity constraints. Therefore, the computation of the main tools needed for the Riemannian optimization procedure may not be obtained as a direct derivation from the formulae in \cite{Reversible}.
\end{remark}

\subsection{Selecting a different pattern from the original matrix}

Analogously to the constrained approach of Section~\ref{subsec:euclidean_optimizer_sub pattern}, we define a Riemannian optimizer for the same setting. Given a reversible matrix $P$ with pattern $\mathcal{P}$, we now consider perturbations $\Delta$ with a (possibly different) pattern $\mathcal{S}$, reflecting situations in which only selected entries of $P$ may be modified.

The objective \eqref{eq:functional_min_X} remains unchanged, but the manifold $\mathcal{M}_{\boldsymbol{\pi}}$ must be adapted to accommodate sub-patterns of $P$. This leads to a feasible set whose entries are partly fixed and must remain so throughout optimization, by imposing an additional constraint on the matrices $D_{\hat{\boldsymbol{\pi}}}(P+\Delta)D_{\hat{\boldsymbol{\pi}}}^{-1}$.

Indeed, the detailed balance condition \eqref{eq:detailed_balance} leads to the construction of the symmetric matrices $X = D_{\boldsymbol{\hat \pi}}( P +\Delta) D_{\boldsymbol{\hat \pi}}^
{-1}$, where $P + \Delta$ is a reversible matrix. In this section, we are interested in modifying selected entries, namely the ones belonging to the pattern $\mathcal{S}$, and leaving the rest unchanged. 

Let us consider a pair $\left\lbrace i,j\right\rbrace \not \in \mathcal{S}$. Then, the corresponding $(P+\Delta)_{ij}$ must be chosen to be equal to $P_{ij}$. In the symmetrized version  $X$ obtained via the scaling by $D_{\boldsymbol{\hat \pi}}$, this can be obtained fixing $X_{ij} = \left[D_{\boldsymbol{\hat \pi}} (P + \Delta ) D_{\boldsymbol{\hat \pi}}^{-1}\right]_{ij}$, that is choosing $X_{ij} = \frac{\hat{\boldsymbol{\pi}}_i}{\hat{\boldsymbol{\pi}}_j} P_{ij}$.

Observe that the symmetry of $X$ is satisfied, since, for $\left\lbrace i,j\right\rbrace \not \in \mathcal{S}$,  we have that $X_{ij} = \frac{\hat{\boldsymbol{\pi}}_i}{\hat{\boldsymbol{\pi}}_j} P_{ij} = \frac{\hat{\boldsymbol{\pi}}_i}{\hat{\boldsymbol{\pi}}_j} \frac{\boldsymbol{\pi}_j}{\boldsymbol{\pi}_i} P_{ji} = \frac{\hat{\boldsymbol{\pi}}_j}{\hat{\boldsymbol{\pi}}_i} P_{ji} = X_{ji}$, where we used again~\eqref{eq:detailed_balance}.

\begin{remark}\label{rem:assumption}
   In this section, we assume that the matrix $P^0$ built as in \eqref{eq:P0} is such that $\| P^0 \mathbf{1}\|_{\infty} <1$; see also Remark~\ref{remark: only one}.
\end{remark}

Therefore, the optimization problem has as minimizing functional $f(X)$ in \eqref{eq:functional_min_X} and as non-empty feasible set (see Remark~\ref{remark: only one}) the {following one}: %
\begin{equation}\label{eq:manifold_fixed_P}
\begin{split}
   \mathcal{M}_{P,\boldsymbol{\pi}} = & \left\lbrace X \in \mathbb{R}^{n\times n}_{\rm{exact}}(\mathcal{P} \cup \mathcal{S}) : X = X^\top, X\hat{\boldsymbol{\pi}} = \hat{\boldsymbol{\pi}}, %
   X_{ij} > 0 \; \mbox{if } \left\lbrace i,j\right\rbrace \in \mathcal{S}, \right.\\ & \left. \quad X_{ij}= \frac{\hat{\boldsymbol{\pi}}_i}{\hat{\boldsymbol{\pi}}_j} P_{ij}   \mbox{ if } \left\lbrace i,j\right\rbrace  \not\in \mathcal{S} \right\rbrace.
\end{split}
\end{equation}

Since the starting matrix $P$ is reversible, we automatically get that $X_{ij} = X_{ji}$ for each $\left\lbrace i,j\right\rbrace  \not\in \mathcal{S}$. Having solved the optimization for $X$, we then recover the reversible stochastic matrix $P+\Delta$ and the perturbation $\Delta$ by inverting the change of variables.

We observe that the set~\eqref{eq:manifold_fixed_P} is a manifold. Firstly, we have that the set
\[
\widetilde{\mathcal{M}} = \left\lbrace X \in \mathbb{R}^{n\times n}_{\rm{exact}}(\mathcal{P} \cup \mathcal{S}) : X = X^\top, X\hat{\boldsymbol{\pi}} = \hat{\boldsymbol{\pi}}, %
   X_{ij}= \frac{\hat{\boldsymbol{\pi}}_i}{\hat{\boldsymbol{\pi}}_j} P_{ij}   \mbox{ if } \left\lbrace i,j\right\rbrace  \not\in \mathcal{S} \right\rbrace
   \]
is an embedded manifold. Indeed, we may construct the local defining function $h$ \cite[Definition 3.10]{BoumalBook}, as follows: 
\[
h(X) = \left( h^{(1)}(X), h^{(2)}(X), h^{(3)}(X) \right),
\]
where $h^{(1)}(X)$ collects the symmetric constraints $X_{ji} - X_{ij}$, for $i<j$, $h^{(2)}(X)$ the ones in the fixed eigenvector $\boldsymbol{\hat \pi}$, that are $\sum_{j} X_{ij} \boldsymbol{\hat \pi}_j - \boldsymbol{\hat \pi}_i$, for $i=1,\ldots,n$, and $h^{(3)}(X)$ fixes $X_{ij} - \frac{\hat{\boldsymbol{\pi}}_i}{\hat{\boldsymbol{\pi}}_j} P_{ij}$, for $\left\lbrace i,j\right\rbrace  \not\in \mathcal{S}$. Secondly, we observe that 
\[
\mathcal{M}_{P,\boldsymbol{\pi}}  = \widetilde{\mathcal{M}} %
\cap 
\left\lbrace \sum_{ \{i,j\} \in \mathcal{S}} \alpha_{i,j}\, \mathbf{e}_i \mathbf{e}_j^\top + \sum_{\{p,q\} \not\in \mathcal{S}} \beta_{p,q}\, \mathbf{e}_p \mathbf{e}_q^\top \,:\, \alpha_{i,j} > 0, \quad \beta_{p,q} \in \mathbb{R} \right\rbrace
\]
is an open set of $\widetilde{\mathcal{M}}$. This implies that $\mathcal{M}_{P,\boldsymbol{\pi}}$ is a manifold.

Similarly to the manifold $\mathcal{M}_{\boldsymbol{\pi}}$ in \eqref{eq:manifold_sparse_reversible}, the manifold $\mathcal{M}_{P,\boldsymbol{\pi}}$ is not a submanifold of the one presented in~\cite{Reversible}, since the entrywise positivity constraint is not satisfied. Hence, in the following, we need to derive the main optimization tools via direct computation, which do not straightforwardly derive from the ones in \cite{Reversible}.

\begin{remark}
    Observe that the set $\mathcal{F}$ defined in \eqref{eq:fixed-feasible} is strictly connected with the manifold $\mathcal{M}_{P,\boldsymbol{ \pi}}$. Indeed, it can be seen that, for matrices ${X} \in \mathbb{R}_{\text{exact}}^{n\times n}(\mathcal{P}\cup \mathcal{S})$, the mapping $\mathcal{L} : {X} \mapsto D_{\hat{\boldsymbol{\pi}}}^{-1} {X} D_{\hat{\boldsymbol{\pi}}}$ transforms elements in $\mathcal{M}_{P,\boldsymbol{\pi}}$ into elements in $\mathcal{F}$. Nevertheless, it is worth noticing that we restrict the feasible set to matrices in $\mathbb{R}_{\text{exact}}^{n\times n}(\mathcal{P}\cup \mathcal{S})$, rather than $\mathbb{R}^{n\times n}(\mathcal{P}\cup \mathcal{S})$. This choice guarantees to work on an open manifold without boundary, which is particularly convenient for the development of Riemannian optimization methods.

    Experimentally, this construction can be mitigated by progressively modifying the set $\mathcal{S}$ of nonzero elements, in the numerical solution of the optimization problem, and adjusting the pattern accordingly. To this end, we propose an experimental strategy described in the subsequent Section~\ref{sec:pattern}.
\end{remark}

\subsection{Main tools for the Riemannian manifold}
\label{sec:main tools}

As anticipated, we may minimize the functional in~\eqref{eq:functional_min_X}, optimizing over a suitable manifold. In particular, we observe that the manifolds $\mathcal{M}_{\boldsymbol{\pi}}$ in \eqref{eq:manifold_sparse_reversible} and $\mathcal{M}_{P,\boldsymbol{\pi}}$ in \eqref{eq:manifold_fixed_P} coincide if $\mathcal{P} \equiv \mathcal{S}$, since this implies $X_{ij} = 0$ for $\left\lbrace i,j\right\rbrace \not\in  \mathcal{S}$. Therefore, we may construct and implement the main tools for the Riemannian optimization procedure by focusing on the manifold $\mathcal{M}_{P,\boldsymbol{\pi}}$ and taking $\mathcal{M}_{\boldsymbol{\pi}}$ as its special case. 

We start the construction with the characterization of the tangent space of the manifold $\mathcal{M}_{P,\boldsymbol{\pi}}$.

\begin{lemma}
The tangent space $\mathcal{T}_X \mathcal{M}_{P,\boldsymbol{\pi}}$ at $X\in \mathcal{M}_{P,\boldsymbol{\pi}}$ is the set
\begin{equation*}
   \mathcal{T}_X \mathcal{M}_{P,\boldsymbol{\pi}} = \left\lbrace \xi_X \in \mathbb{R}^{n\times n}(\mathcal{S}) : \xi_X = \xi_X^\top, \xi_X \hat{\boldsymbol{\pi}}=0 %
   \right\rbrace.
\end{equation*}   
\end{lemma}

\begin{proof}
   Consider a smooth curve $X(t) \in \mathcal{M}_{P,\boldsymbol{\pi}}$, with $X(0)=X$, for sufficiently small values of $t$. Then, we have
   \[
   \begin{array}{c}
      X(t) = X(t)^\top \implies \dot{X}(t) = \dot{X}(t)^{\top}, \;       X(t)\hat{\boldsymbol{\pi}} = \hat{\boldsymbol{\pi}} \implies \dot{X}(t)\hat{\boldsymbol{\pi}} = 0, \\
       X_{ij}(t) = \frac{\hat{\boldsymbol{\pi}}_i}{\hat{\boldsymbol{\pi}}_j} P_{ij}, \;\forall t, \forall \left\lbrace i,j\right\rbrace  \not\in \mathcal{S} \implies \dot{X}_{ij}(t)= 0 ,\; \forall \left\lbrace i,j\right\rbrace  \not\in \mathcal{S}.
   \end{array}
   \]
   This implies that $\mathcal{T}_{X}\mathcal{M}_{P,\boldsymbol{\pi}} \subseteq \left\lbrace \xi_X \in \mathbb{R}^{n\times n}(\mathcal{S})  : \xi_X = \xi_X^\top, \xi_X \hat{\boldsymbol{\pi}}=0 %
   \right\rbrace$.

   On the other hand, given $W$ in the set at the right hand side, the curve $\gamma(t) = X + tW$ belongs to $\mathcal{M}_{P,\boldsymbol{\pi}}$ for sufficiently small values of $t$. Since $\gamma(0)=X$ and $\gamma'(0)=W$, we have the inclusion $\mathcal{T}_{X}\mathcal{M}_{P,\boldsymbol{\pi}} \supseteq \left\lbrace \xi_X \in \mathbb{R}^{n\times n}(\mathcal{S}) : \xi_X = \xi_X^\top, \xi_X \hat{\boldsymbol{\pi}}=0 %
   \right\rbrace$.
\end{proof}

Observe that the tangent space only takes into account the set of admissible directions $X(t)$, for sufficiently small values of $t$. Therefore, since the entries $X_{ij}(t)$ are constantly equal for every $t$ and for $\left\lbrace i,j\right\rbrace \not \in \mathcal{S}$, the first derivative $\dot{X}_{ij}(t)$ does not reflect any changes outside the set of indices $\mathcal{S}$. Moreover, using this result, we get that the tangent space $\mathcal{T}_{X} \mathcal{M}_{P,\boldsymbol{\pi}}$ does not depend on the point $X$, but only on the chosen pattern and the stationary distribution.

As a sanity check, we observe that the tangent space $\mathcal{T}_{X} \mathcal{M}_{P,\boldsymbol{\pi}}$ is equal to the tangent space $\mathcal{T}_X \mathcal{M}_{\boldsymbol{\pi}}$. Then, we can endow it with the same inner product employed for the manifold \eqref{eq:manifold_sparse_reversible}, that is: given $\xi_X, \eta_X \in \mathcal{T}_{X} \mathcal{M}_{P,\boldsymbol{\pi}}$
\[
\left\langle \xi_X, \eta_X \right\rangle_X = \sum_{X_{ij}\neq 0} \frac{(\xi_X)_{ij}(\eta_X)_{ij}}{X_{ij}}.
\]
The characterization of the tangent space can be used in constructing the projection $\Pi_X: \mathbb{R}^{n\times n}\mapsto \mathcal{T}_X\mathcal{M}_{P,\boldsymbol{\pi}}$, by characterizing $\mathcal{T}_X\mathcal{M}_{P,\boldsymbol{\pi}}^{\perp}$, the orthogonal complement of the tangent space in $\mathbb{R}^{n\times n}(\mathcal{S})$. %
In this setting, we may prove the following.
\begin{lemma}
The orthogonal complement $\mathcal{T}_X ^{\perp} \mathcal{M}_{P,\boldsymbol{\pi}}$ to the tangent space $\mathcal{T}_X \mathcal{M}_{P,\boldsymbol{\pi}}$ can be expressed as
   \begin{equation}
   \label{eq:orthogonal_complement}
   \mathcal{T}_X ^{\perp} \mathcal{M}_{P,\boldsymbol{\pi}} = \left\lbrace \xi_X^{\perp} \in \mathbb{R}^{n\times n}: \xi_X^{\perp} = \left(\boldsymbol{\alpha}\boldsymbol{\hat\pi}^{\top} + \boldsymbol{\hat \pi} \boldsymbol{\alpha}^{\top} \right)\odot (X \odot S) \right\rbrace,
\end{equation}
where $\boldsymbol{\alpha} \in \mathbb{R}^{n}$ and the matrix $S\in \mathbb{R}^{n\times n}$ is defined as
\begin{equation}\label{eq:sparse_matrix}
   S_{ij}= \left\lbrace \begin{array}{c}
        1, \qquad \mbox{for } \left\lbrace i,j\right\rbrace \in \mathcal{S},\\
        0, \qquad \mbox{for } \left\lbrace i,j\right\rbrace \not\in \mathcal{S}.
   \end{array} \right.
\end{equation}
\end{lemma}
\begin{proof}
   Firstly, we prove the inclusion $ \supseteq$. %
   Given $\xi_X \in \mathcal{T}_{X}\mathcal{M}_{P,\boldsymbol{\pi}}$ and $Z=(\boldsymbol{\alpha}\boldsymbol{\hat \pi}^{\top} + \boldsymbol{\hat \pi} \boldsymbol{\alpha}^{\top})\odot (X \odot S)$, we get
\begin{align*}
       \left\langle \xi_X, Z \right\rangle_X &=  \sum_{X_{ij}\neq 0} \frac{(\xi_X)_{ij}{\left((\boldsymbol{\alpha}\boldsymbol{\hat \pi}^{\top} + \boldsymbol{\hat \pi} \boldsymbol{\alpha}^{\top})\odot (X \odot S)\right)}_{ij}}{X_{ij}} \\
       =& \sum_{X_{ij}\neq 0} \frac{(\xi_X)_{ij} \boldsymbol{\alpha}_i \boldsymbol{\hat \pi}_j X_{ij} S_{ij}}{X_{ij}} + \sum_{X_{ij}\neq 0} \frac{(\xi_X)_{ij} \boldsymbol{\hat \pi}_i \boldsymbol{\alpha}_j X_{ij} S_{ij}}{X_{ij}} \\
       & = \boldsymbol{\alpha}^{\top} (\xi_X \odot S) \boldsymbol{\hat \pi} + \boldsymbol{\hat \pi}^{\top} (\xi_X \odot S) \boldsymbol{\alpha} = \boldsymbol{\alpha}^{\top} \underbrace{\xi_X \boldsymbol{\hat \pi}}_{=0} +  \underbrace{\boldsymbol{\hat \pi}^{\top} \xi_X}_{=0} \boldsymbol{\alpha} =0,
   \end{align*}
   where we used the fact that $\xi_X$ is symmetric and has pattern $\mathcal{S}$. Moreover, we observe that any matrix $ \left(\boldsymbol{\alpha}\boldsymbol{\hat\pi}^{\top} + \boldsymbol{\hat \pi} \boldsymbol{\alpha}^{\top} \right)\odot (X \odot S)$ is symmetric and has pattern $\mathcal{S}$, and then it is contained in $\mathbb{R}^{n\times n}(\mathcal{S})$. The inclusion holds since the two sets have the same dimension, as vector spaces. Indeed, the set at the right hand side in~\eqref{eq:orthogonal_complement} is a vector space of dimension $\left| \mathcal{S}\right| -n$. The dimension of the orthogonal complement of $\mathcal{T}_{X}\mathcal{M}_{P,\boldsymbol{\pi}}$ in $\mathbb{R}^{n\times n}(\mathcal{S})$ is given by $\text{dim}(\mathbb{R}^{n\times  n}(\mathcal{S})) - \text{dim}(\mathcal{T}_{X} \mathcal{M}_{P,\boldsymbol{\pi}}) = \left|  \mathcal{S} \right| -n$.
\end{proof}

The expression for the elements in $\mathcal{T}_X \mathcal{M}_{P,\boldsymbol{\pi}}$ leads to the possibility to define a projection from the embedding space $\mathbb{R}^{n\times n}$ to the tangent space of the manifold. This step will be employed in the construction of the Riemannian gradient of a function defined on $\mathcal{M}_{P,\boldsymbol{\pi}}$, given its Euclidean counterpart.

The projection $\Pi_X$ from the embedding space to the tangent space can be written~as
\begin{align*}
   \Pi_X : \mathbb{R}^{n\times n} \mapsto \mathcal{T}_X \mathcal{M}_{P,\boldsymbol{\pi}} \text{ with }
   Z \mapsto \frac{(Z+Z^{\top})\odot S}{2} - (\boldsymbol{\alpha}\boldsymbol{\hat \pi}^{\top} + \boldsymbol{\hat \pi} \boldsymbol{\alpha}^{\top})\odot (X \odot S),
\end{align*}
where $S$ is defined as in \eqref{eq:sparse_matrix}
and the vector $\boldsymbol{\alpha}$ can be derived as the solution of 
\[
\frac{(Z+Z^\top)\odot S}{2}\boldsymbol{\hat \pi} = \left( (\boldsymbol{\alpha}\boldsymbol{\hat \pi}^{\top} + \boldsymbol{\hat \pi} \boldsymbol{\alpha}^{\top}) \odot (X \odot S)\right) \boldsymbol{\hat \pi}.
\]
Manipulating the previous expression, using computations similar to the one proposed in \cite[Theorem 8]{Douik} and \cite[Lemma 3.3]{Reversible}, we get that the vector $\boldsymbol{\alpha}$ is the solution of the linear system
\begin{equation}
\label{eq:linear_system_proj}
   \frac{(Z+Z^\top)\odot S}{2}\boldsymbol{\hat \pi} = 
   \left( \operatorname{diag}((X \odot S) \boldsymbol{\pi}) + D_{\boldsymbol{\hat \pi}}(X \odot S) D_{\boldsymbol{\hat \pi}})\right)\boldsymbol{\alpha}.
\end{equation}

\begin{proposition}
   The matrix $\left( {\operatorname{diag}}((X \odot S) \boldsymbol{\pi}) + D_{\boldsymbol{\hat \pi}}(X \odot S) D_{\boldsymbol{\hat \pi}})\right)$ in~\eqref{eq:linear_system_proj} is invertible for any $X \in \mathcal{M}_{P,\boldsymbol{\pi}}$, and $S$ as in~\eqref{eq:sparse_matrix}.
\end{proposition}

\begin{proof}
   Assume by contradiction that there exists a nonzero vector $\mathbf{v}\in \mathbb{R}^{n}$ such that
\begin{align*}
   \operatorname{diag}((X \odot S) \boldsymbol{\pi})\mathbf{v} + D_{\boldsymbol{\hat \pi}}(X \odot S) D_{\boldsymbol{\hat \pi}} \mathbf{v}= 0 \iff  \left( \operatorname{diag}((X \odot S) \boldsymbol{\pi}) \right)^{-1} D_{\boldsymbol{\hat \pi}}(X \odot S)D_{\boldsymbol{\hat \pi}} \mathbf{v}=  - \mathbf{v},
\end{align*}
by observing that the matrix $\operatorname{diag}((X \odot S) \boldsymbol{\pi})$ is invertible since $(X \odot S)_{ii}\neq 0$, $(X \odot S) \geq 0$, and $\boldsymbol{\pi} > 0$.

As an intermediate step, we observe that 
\begin{align}
\label{eq:collatz}
   \| \left( \operatorname{diag}(X \boldsymbol{\pi}) \right)^{-1} D_{\boldsymbol{\hat \pi}}X D_{\boldsymbol{\hat \pi}} \|_{\infty} = \max_{1\leq i \leq n} \frac{\boldsymbol{\pi}_i}{\left( X \boldsymbol{\pi} \right)_i} = \min_{1\leq i \leq n} \frac{\left( X \boldsymbol{\pi} \right)_i}{\boldsymbol{\pi}_i} \leq \rho (X) = 1,
\end{align}
where the first equality follows by the nonnegativity of the matrix, and by noting that $\left( \operatorname{diag}(X \boldsymbol{\pi}) \right)^{-1} D_{\boldsymbol{\hat \pi}}X D_{\boldsymbol{\hat \pi}} \mathbf{1} = \left( \operatorname{diag}(X \boldsymbol{\pi}) \right)^{-1} \boldsymbol{\pi}$; while the last inequality follows from the Collatz--Wielandt formula \cite[\S 8]{Meyer}. 

From~\eqref{eq:collatz}, we get that $X\boldsymbol{\pi} \geq \boldsymbol{\pi}$ entrywise. Since $X$ irreducible and nonnegative, with $\rho (X) =1$, from \cite[Example 8.3.1]{Meyer} we get that $X \boldsymbol{\pi} = \boldsymbol{\pi}$, which is a contradiction since $\boldsymbol{\hat \pi}$ is the unique eigenvector associated with $\rho(X)$ and $\boldsymbol{\hat \pi} > \boldsymbol{\pi}$ entrywise, hence it cannot be a multiple of $\boldsymbol{\pi}$ due to the fact that $\boldsymbol{\hat \pi} = \boldsymbol{\pi}^{\nicefrac{1}{2}}$.

Moreover, since the infinity norm of the matrix $\left( \operatorname{diag}((X \odot S) \boldsymbol{\pi}) \right)^{-1} D_{\boldsymbol{\hat \pi}} (X \odot S) D_{\boldsymbol{\hat \pi}}$ can be bounded from above as
\begin{align*}
   \| \left( \operatorname{diag}((X \odot S)\boldsymbol{\pi}) \right)^{-1} D_{\boldsymbol{\hat \pi}}(X \odot S) D_{\boldsymbol{\hat \pi}} \|_{\infty} &= \max_{1\leq i \leq n} \frac{(D_{\boldsymbol{\hat \pi}}(X \odot S) \boldsymbol{\hat\pi})_i}{((X \odot S)\boldsymbol{\pi})_i} \leq\max_{1\leq i \leq n} \frac{(D_{\boldsymbol{\hat \pi}}X {\boldsymbol{\hat \pi}})_i }{((X \odot S)\boldsymbol{\pi})_i} \\
   &= \min_{1\leq i \leq n} \frac{((X \odot S)\boldsymbol{\pi})_i}{(D_{\boldsymbol{\hat \pi}}X {\boldsymbol{\hat \pi}})_i} \leq \min_{1\leq i \leq n} \frac{\left( X \boldsymbol{\pi} \right)_i}{(D_{\boldsymbol{\hat \pi}}X {\boldsymbol{\hat \pi}})_i}\\
   &= \| \left( \operatorname{diag}(X \boldsymbol{\pi}) \right)^{-1} D_{\boldsymbol{\hat \pi}}X D_{\boldsymbol{\hat \pi}} \|_{\infty},
\end{align*}
we get that $\rho\left( \left( \operatorname{diag}((X \odot S) \boldsymbol{\pi}) \right)^{-1} D_{\boldsymbol{\hat \pi}} (X \odot S) D_{\boldsymbol{\hat \pi}} \right)<1$, which concludes the proof.
\end{proof}

In order to employ first order Riemannian optimizers, such as the conjugate gradient~\cite{MR3325229} or Barzilai-Borwein method~\cite{Iannazzo}, we need the computation of the Riemannian gradient~\cite[Deﬁnition 3.58]{BoumalBook} of the function $f(X)$ in \eqref{eq:functional_min_X}. Denote by $\operatorname{grad} f$ and $\operatorname{Grad} f$ the Riemannian and Euclidean gradient~\cite[eq. (3.19)]{BoumalBook} of the function $f:\mathcal{M}_{P,\boldsymbol{\pi}} \mapsto \mathbb{R}$, respectively. By definition of Euclidean gradient, we get that 
\[
  \left\langle \operatorname{Grad} f(X), \xi \right\rangle = \mathrm{D} f(X) [\xi], \quad \forall \xi \in \mathbb{R}^{n\times n},
\]
where $\left\langle \cdot , \cdot \right\rangle$ denotes the Frobenius inner product. The definition of the Riemannian metric in \eqref{eq:fisher_sparse} and the formula for the orthogonal complement of the tangent space $\mathcal{T}_{X} \mathcal{M}_{P,\boldsymbol{\pi}}$ lead to the expression
\[
\left\langle \operatorname{Grad} f(X), \xi_X \right\rangle =   \left\langle {\Pi}_X (\operatorname{Grad} f(X) \odot X),  \xi_X \right\rangle_X = \mathrm{D} f(X) [\xi_X],
\]
for each $\xi_X \in \mathcal{T}_X{\mathcal{M}_{P,\boldsymbol{\pi}}}$. From which we get that the Riemannian gradient $\operatorname{grad} f(X)$ is equal to ${\Pi}_X (\operatorname{Grad} f(X) \odot X)$.

The last tool we need for employing first order methods on a manifold is the possibility to perform a retraction from the tangent space to the manifold~\cite[Definition~3.47]{BoumalBook}. In the following, we construct a first-order retraction of a point $\xi_X$ in the tangent space onto the manifold $\mathcal{M}_{P,\boldsymbol{\pi}}$. As an intermediate step, we need the following lemma, which describes the possibility to balance a symmetric matrix. The proof relies on a 
modified version of Sinkorn--Knopp~\cite{MR210731,MR161868,Reversible}, and the condition in~\cite[Theorem~2.1]{MR2399579}. 

\begin{lemma}\label{lem:modified_SK}
   Let $P \in\mathbb{R}^{n\times n}$ be a reversible stochastic matrix with stationary distribution $\boldsymbol{\pi}$, and $\mathcal{S}$ a pattern for which $P^0$ in~\eqref{eq:P0} is such that $\|P^0\mathbf{1}\|_{\infty} < 1$. Then, for any nonnegative symmetric matrix $A\in \mathbb{R}^{n \times n}(\mathcal{S})$, $A \neq 0$ with total support---i.e., such that all its nonzero elements lie on a {diagonal with positive elements} %
   ---there exists a diagonal matrix $D$ with positive diagonal entries such that
   $
   DAD \boldsymbol{\hat \pi} = \boldsymbol{\hat \pi} - \boldsymbol{\beta}$, with  $\boldsymbol{\beta}=P^0 \hat{\boldsymbol{\pi}}.
   $
\end{lemma}
\begin{proof}
   Note that the vector at the right hand side has entries different from zero. %
   Then, applying a modified version of the Sinkorn--Knopp to the matrix $\hat A = D_{\boldsymbol{\hat \pi}} AD_{\boldsymbol{\hat \pi}}$, as in \cite[Lemma~3.7]{Reversible}, we obtain two diagonal matrices $D_1,D_2$ with positive elements on the diagonals such that
\begin{equation}
\label{eq:SK_relation}
D_1 \hat A D_2 \mathbf{1} = D_{\boldsymbol{\hat \pi}}(\boldsymbol{\hat \pi}- \boldsymbol{\beta}) \; \mbox{and } \mathbf{1}^{\top} D_1 \hat A D_2 = (\boldsymbol{\hat \pi} - \boldsymbol{\beta})^{\top} D_{\boldsymbol{\hat \pi}}.
\end{equation}
From the previous relation, using that $\hat A =D_{\boldsymbol{\hat \pi}} AD_{\boldsymbol{\hat \pi}}$, we get
\[
D_1 A D_2 \boldsymbol{\hat \pi} = \boldsymbol{\hat \pi} - \boldsymbol{\beta}\; \mbox{and } \boldsymbol{\hat \pi}^{\top} D_1 A D_2 = (\boldsymbol{\hat \pi} - \boldsymbol{\beta})^{\top}.
\]
Since $A$ is symmetric, taking the transpose in \eqref{eq:SK_relation}, it also holds
\[
D_2 A D_1 \boldsymbol{\hat \pi} = \boldsymbol{\hat \pi} - \boldsymbol{\beta}\; \mbox{and } \boldsymbol{\hat \pi}^{\top} D_2 A D_1 = (\boldsymbol{\hat \pi} - \boldsymbol{\beta})^{\top}.
\]
Then, we conclude the proof choosing $D = \frac{1}{2}(D_1+D_2)$. %
\end{proof}

Then, given a nonnegative and symmetric matrix $A$, with a prescribed pattern $\mathcal{S}$, we find a diagonal scaling $D$ such that
the matrix $\widetilde P +DAD$ belongs to the manifold $\mathcal{M}_{P,\boldsymbol{\pi}}$, where $\widetilde P_{ij}=P_{ij}\frac{\boldsymbol{\hat\pi}_i}{\boldsymbol{\hat\pi}_j}$, for $\left\lbrace i,j \right\rbrace \not \in \mathcal{S}$ and $\widetilde P_{ij}=0$ elsewhere. Indeed, we may check that
$
\widetilde P \boldsymbol{\hat \pi} + DAD \boldsymbol{\hat \pi} = \boldsymbol{\beta} + \boldsymbol{\hat \pi} - \boldsymbol{\beta} = \boldsymbol{\hat \pi},
$
and consequently, $\boldsymbol{\hat \pi}^{\top} ( \widetilde P + DAD) = \boldsymbol{\hat \pi}^{\top}$, since $\widetilde P$ is symmetric.

\begin{remark}
   Observe that when dealing with  matrices in $\mathcal{M}_{\boldsymbol{\pi}}$, where the sparsity pattern is the one inherited by the original matrix, we have $\boldsymbol{\beta} = 0$ in the previous Lemma.
\end{remark}

A first-order retraction $R_X: \mathcal{T}_X \mathcal{M}_{P,\boldsymbol{\pi}} \mapsto \mathcal{M}_{P,\boldsymbol{\pi}}$ onto the manifold can be constructed via the result \cite[Proposition 4.1.2]{AbsilBook}. In this framework, this means constructing a diffeomorphism $\phi$ such that
\begin{align*}
   \phi : \mathcal{M}_{P,\boldsymbol{\pi}} \times \mathbb{R}^{n}_{>} \mapsto \overline{\mathcal{S}}_P, \text{ with } 
   (\tilde{P}+S,\mathbf{d}) \mapsto \tilde{P} + D_{\mathbf{d}} SD_{\mathbf{d}},
\end{align*}
where $\mathbb{R}_{>}^{n}= \left\lbrace \boldsymbol{x} \in \mathbb{R}^n : \boldsymbol{x}>0\right\rbrace$ and $\overline{\mathcal{S}}_P=\left\lbrace \tilde{P}+ S: S \in \mathbb{R}^{n\times n}_{\rm{exact}}(\mathcal{S}), S_{ij} >0 \text{ for }\left\lbrace i,j \right\rbrace \in \mathcal{S}\right\rbrace$. %
Recall that each element of the manifold $\mathcal{M}_{P,\boldsymbol{\pi}}$ can be written as a sum of the matrix with elements $\tilde{P}_{ij}=P_{ij}\frac{\boldsymbol{\hat \pi}_i}{\boldsymbol{\hat \pi}_j}$ and a matrix $S\geq 0$ in $\mathbb{R}^{n\times n}_{\rm exact}(\mathcal{S})$. Note that $\operatorname{dim}(\mathcal{M}_{P,\boldsymbol{\pi}}) = \operatorname{dim}(\overline{\mathcal{S}}_P) - n$, moreover $\mathbb{R}^n_{>}$ and $\overline{\mathcal{S}}_P$ are manifold. 

The function $\phi$ is a diffeomorphism. Indeed, observe that its inverse $\phi^{-1}$ can be derived employing Lemma~\ref{lem:modified_SK} and the subsequent remark. Then, a retraction $R_X : \mathcal{T}_X\mathcal{M}_{P,\boldsymbol{\pi}} \mapsto \mathcal{M}_{P,\boldsymbol{\pi}}$ can be constructed as $R_X(\xi_X)=\tau_1\left( \phi^{-1}(X+\xi_X)\right)$, where $\tau_1: \mathcal{M}_{P,\boldsymbol{\pi}} \times \mathbb{R}_{>}^n \mapsto \overline{\mathcal{S}}_P$ is the projection onto the first component. Observe that since $X=\tilde{P} + S$, we have that the first order retraction $R_X$ is the identity, for $X+\xi_X>0$. Indeed, we have that
$
(\tilde{P} + S +\, \xi_X)\boldsymbol{\hat \pi} = (X+\xi_X) \boldsymbol{\hat \pi} = \boldsymbol{\hat \pi}, \mbox{ and } (\tilde{P} + S + \xi_X)= (X + \xi_X) = (X  +\xi_X)^{\top}.
$

In order to employ a Riemannian optimization scheme for the functional $f(X)$ in \eqref{eq:functional_min_X}, we need as an additional step the computation of the Euclidean gradient of $f(X)$, which will then be projected in order to find its Riemannian counterpart. In the following, we denote by $f$ the smooth extension of $f$ defined on $\mathcal{M}_{P,\boldsymbol{\pi}}$ to the embedded space $\mathbb{R}^{n\times n}$.

\begin{lemma}
   Let $f:\mathbb{R}^{n \times n} \mapsto \mathbb{R}$ be in the form
   \[
   f(X)= \operatorname{tr} \left( (I - X + {\boldsymbol{\hat \pi}}{\boldsymbol{\hat \pi}}^{\top} )^{-1}\right) + \frac{1}{2}\|D_{\boldsymbol{\hat\pi}}^{-1} X D_{\boldsymbol{\hat \pi}} - P \|_F^2,
   \]
   then the Euclidean gradient $\operatorname{Grad} f(X)$ is
   \[
 \operatorname{Grad} f(X)=  D_{\boldsymbol{\pi}}^{-1} X D_{\boldsymbol{\pi}} - D_{\boldsymbol{\hat \pi}}^{-1} P D_{\boldsymbol{\hat \pi}} + \left((I-X+{\boldsymbol{\hat \pi}}{\boldsymbol{\hat \pi}}^{\top})^2\right)^{-\top}.
   \]
\end{lemma}

\begin{proof}
    The derivative for the term $\|D_{\boldsymbol{\hat\pi}}^{-1} X D_{\boldsymbol{\hat \pi}} - P \|_F^2$ can be found in \cite[Proposition~3.9]{Reversible}. For the second term, we observe that $\operatorname{tr} \left( (I - X + {\boldsymbol{\hat \pi}}{\boldsymbol{\hat \pi}}^{\top} )^{-1}\right) = h(g(X))$, with  $h(Y)= \operatorname{tr}(Y)$ and $g(X)= (I-X+\boldsymbol{\hat \pi} \boldsymbol{\hat \pi}^{\top})^{-1}$, from which we may conclude that
 \begin{align*}
   \mathrm{D} \operatorname{tr} \left( (I - X + {\boldsymbol{\hat \pi}}{\boldsymbol{\hat \pi}}^{\top} )^{-1}\right) [V] %
   =\left\langle V,  \left((I-X+{\boldsymbol{\hat \pi}}{\boldsymbol{\hat \pi}}^{\top})^2\right)^{-\top}\right\rangle, \text{ for } V \in \mathbb{R}^{n\times n}.
   \end{align*}
Then, we conclude that:  
\begin{align*}
       \mathrm{D} f(X)[V] = \left\langle V, D_{\boldsymbol{\pi}}^{-1} X D_{\boldsymbol{\pi}} - D_{\boldsymbol{\hat \pi}}^{-1} P D_{\boldsymbol{\hat \pi}} + \left((I-X+{\boldsymbol{\hat \pi}}{\boldsymbol{\hat \pi}}^{\top})^2\right)^{-\top} \right\rangle,
   \end{align*}
and the gradient is $\operatorname{Grad} f(X) =  D_{\boldsymbol{\pi}}^{-1} X D_{\boldsymbol{\pi}} - D_{\boldsymbol{\hat \pi}}^{-1} P D_{\boldsymbol{\hat \pi}} + \left((I-X+{\boldsymbol{\hat \pi}}{\boldsymbol{\hat \pi}}^{\top})^2\right)^{-\top}$.
\end{proof}

\subsection{Adaptive choice of the pattern}
\label{sec:pattern}
The Riemannian construction we have described is based on the modified Fisher metric in~\eqref{eq:fisher_sparse}. Unlike the classical case, which requires all matrix entries to be nonzero, the proposed modification permits us to designate a set of entries that remain unaffected by the optimization and can therefore be fixed to zero---subject to the conditions stated in Propositions~\ref{pro:pattern-beatrice} and \ref{pro:metropolis-mod}. At first glance, this flexibility might appear sufficient for the task of reducing Kemeny’s constant: enhancing the communicability and reachability among the chain states (i.e., the nodes of the underlying network) seemingly requires only the addition or reinforcement of connections. However, this intuition is misleading because of the so-called \emph{Braess paradox}~\cite{BraessOriginal}; see the example in Fig.~\ref{fig:braess_paradox}.
\begin{figure}[h]
   \centering
   \includegraphics[width=0.66\columnwidth]{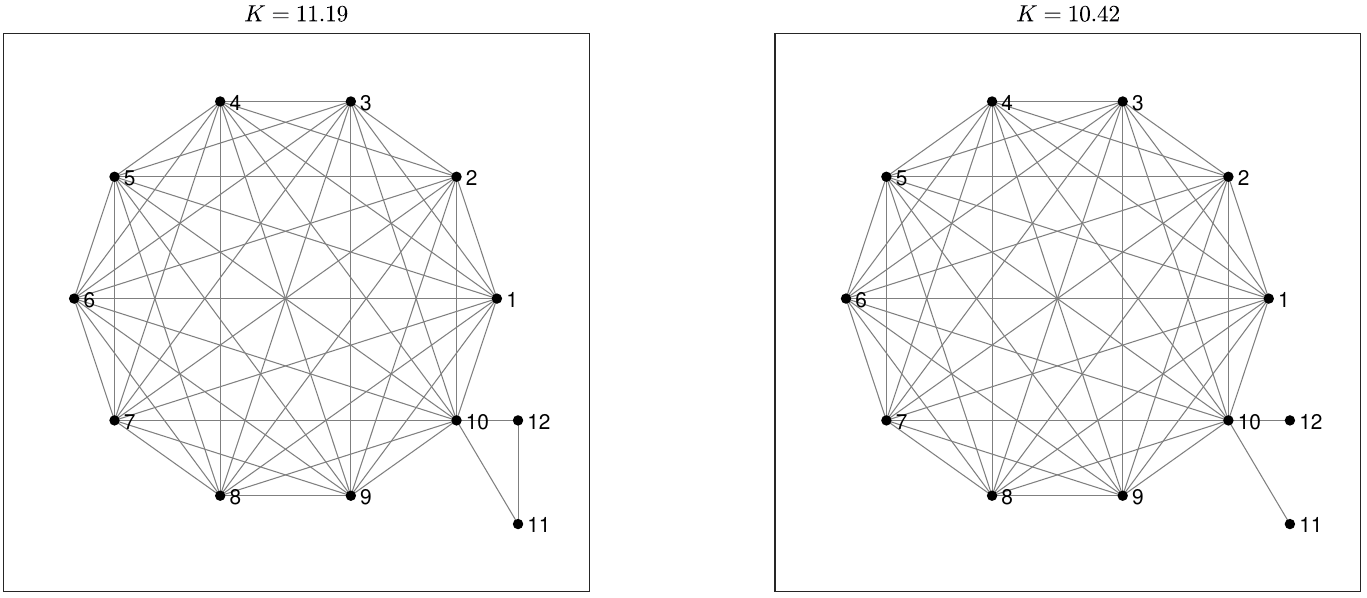}
   \caption{Example of the Braess paradox, for the Markov chains associated with the two graphs.
   In the right panel, the edge between nodes $11$ and $12$ is removed, which decreases Kemeny’s constant from $11.19$ to $10.42$. This illustrates an instance of \emph{pendent twins}, as described in~\cite{Ciardo}.} %
   \label{fig:braess_paradox}
\end{figure}

The Braess paradox~\cite{BraessOriginal,Piqou20171} describes the counterintuitive phenomenon whereby adding an edge or increasing the capacity of a network can worsen overall performance. In the setting of Markov chains, a similar effect can occur: adding or strengthening transitions does not necessarily decrease Kemeny’s constant, and in fact can lead to its increase. This arises because network performance is not determined solely by local improvements in connectivity, but by the global balance of flows and transition probabilities. 

The Riemannian formulation suggests solving \eqref{eq:optimization-fixed-pi-fixed} directly over the manifold $\mathcal{M}_{P,\boldsymbol{\pi}}$, where the optimizer is unique. Computationally, edges associated with the Braess paradox naturally shrink to zero; to prevent numerical issues in evaluating~\eqref{eq:fisher_sparse} when entries fall below machine precision, we adopt an adaptive pattern-selection strategy.
Given a stochastic matrix $P$ with stationary distribution $\boldsymbol{\pi}$ and pattern $\mathcal{S}$, we first optimize over $\mathcal{M}_{P,\boldsymbol{\pi}}$ using the full pattern. Whenever an entry of the current iterate $X$ drops below machine precision $u$, we remove its location from $\mathcal{S}$ and continue optimization on the reduced pattern $\mathcal{S} \setminus \{\{i,j\} \mid X_{ij} \le u,\; i \neq j\}$. Since the minimizer of \eqref{eq:optimization-fixed-pi-fixed} is unique, this pruning cannot introduce spurious solutions. The full procedure is summarized in Algorithm~\ref{alg:adaptive_riemm}.

\begin{algorithm}[htb]
	\caption{Adaptive selection of the pattern}
	\label{alg:adaptive_riemm}
	\begin{algorithmic}
       \State \textbf{Input:} Reversible stochastic matrix $P$, stationary distribution $\boldsymbol{\pi}$ of $P$, prescribed pattern $\mathcal{S}$ %
       \State \textbf{Output:} Reversible stochastic matrix $\hat{P}$ which minimizes the functional \eqref{eq:functional_min_X}
       \State Construct the manifold $\mathcal{M}_{P,\boldsymbol{\pi}}$ as in~\eqref{eq:manifold_fixed_P}
       \State Set tolerance on the norm of the gradient \texttt{tol}= $10^{-3}$
       \While{\texttt{tol}$ > 10^{-12}$}
           \State Compute solution $X$ to problem~\eqref{eq:functional_min_X}, %
           with inner tolerance \texttt{tol}
           \State Choose the pattern $\mathcal{S} \gets \mathcal{S} \; \setminus \; \bigl\{\,\{i,j\} \;\bigm|\; X_{ij} \leq u \land i \neq j \bigr\}$
           \State Choose $X$ as starting guess for the optimization
           \State Decrease the tolerance \texttt{tol} $\gets$ $10^{-3}\cdot\texttt{tol}$
       \EndWhile
       \State Compute the solution as $\hat{P} = D_{\boldsymbol{\hat{\pi}}}^{-1} X D_{\boldsymbol{\hat{\pi}}}$ %
	\end{algorithmic}
\end{algorithm}
This approach mirrors optimization on fixed-rank matrix manifolds, where the chosen rank may initially be suboptimal and is adjusted adaptively, as in~\cite{MR4359465}. Here, varying rank is replaced by varying sparsity patterns: instead of moving between manifolds of different rank, we traverse manifolds of the form~\eqref{eq:manifold_fixed_P}, each defined by a distinct nonzero structure. For solving \eqref{eq:functional_min_X}, any first-order Riemannian method may be used. In Section~\ref{sec:numerical-examples}, we employ the Riemannian Conjugate Gradient (CG) method~\cite{MR3325229} and the Barzilai--Borwein (BB) line-search Riemannian Gradient method~\cite{Iannazzo}.

As an illustrative example, we consider a $10 \times 10$ with a pattern $\mathcal{S}$, with cardinality $54$. Running the experiment in Algorithm~\ref{alg:adaptive_riemm}, we remove from the pattern a few elements, as shown in Fig.~\ref{fig:patter_removal}.

The \ref{fig:original_pattern} panel reports the initial pattern, while the \ref{fig:removed_tol1} and \ref{fig:removed_tol2} panels report the elements that were adaptively selected for elimination from the pattern by the Algorithm~\ref{alg:adaptive_riemm} for two different tolerance values of $10^{-3}$ and $10^{-6}$.

While the idea underlying the adaptive pattern selection may appear straightforward, its introduction was useful to mitigate the numerical errors arising from entries of small magnitude. Avoiding working with such entries is always desirable, in this setting, this concern becomes particularly critical in the optimizer steps that require the retraction $R_X$, described in Section~\ref{sec:main tools}. Specifically, in the implementation of the retraction, we adopt a slightly different first-order retraction, that is
\begin{equation*}
    R_X(\xi_X) = \mathcal{P} \left( X \odot \mathrm{exp}(\xi_X \odiv X) \right), \quad \forall \xi_X \in \mathcal{T}_X \mathcal{M}_{P,\boldsymbol{ \pi}},
\end{equation*}
where the exponential $\exp(\cdot)$ is performed entrywise and $\mathcal{P}: {\mathbb{R}^{n\times n}(\mathcal{S})} \to \mathcal{M}_{P,\boldsymbol{ \pi}}$ represents the application of the modified Sinkhorn theorem proposed in Lemma \ref{lem:modified_SK}. This modification does not alter the theoretical results, while improving the quality of the numerical output, as previously observed in~\cite{Douik, Reversible, DurastanteMeini}. 

In the illustrative example presented in Fig.~\ref{fig:patter_removal}, the removed entries (denoted as  \textcolor{red}{$\blacksquare$}  and \textcolor{green}{$\blacksquare$}) have a magnitude of $10^{-30}$, making
any entry-wise division by them numerically unsafe. Then, adapting the pattern obtained by removing such entries enhances the robustness of the algorithm, without introducing spurious solutions.

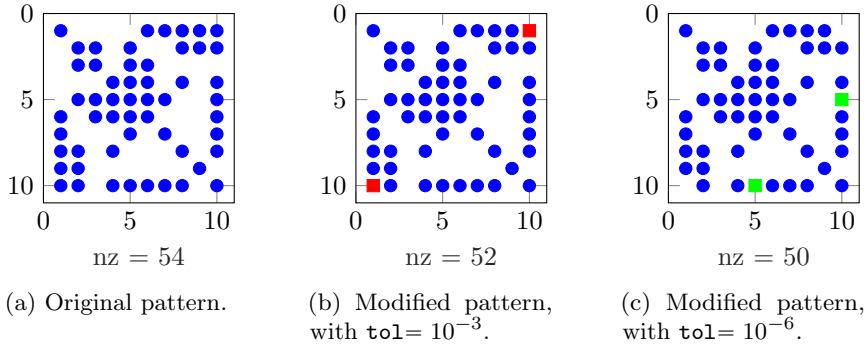
\begin{figure}[hb]
\label{fig:pattern_removal}
   \centering
   \begin{subfigure}[t]{0.24\textwidth}
       \centering
\begin{tikzpicture}

\begin{axis}[%
width=0.8\textwidth,
height=0.8\textwidth,
at={(1.236in,0.481in)},
scale only axis,
xmin=0,
xmax=11,
xlabel style={font=\color{white!15!black}},
xlabel={nz = 54},
y dir=reverse,
ymin=0,
ymax=11,
axis background/.style={fill=white}
]
\addplot [color=blue, only marks, mark size=2.3pt, mark=*, mark options={solid, blue}, forget plot]
 table[row sep=crcr]{%
1	1\\
1	6\\
1	7\\
1	8\\
1	9\\
1	10\\
2	2\\
2	3\\
2	5\\
2	8\\
2	9\\
2	10\\
3	2\\
3	3\\
3	5\\
3	6\\
4	4\\
4	5\\
4	6\\
4	8\\
4	10\\
5	2\\
5	3\\
5	4\\
5	5\\
5	6\\
5	7\\
5	10\\
6	1\\
6	3\\
6	4\\
6	5\\
6	6\\
6	10\\
7	1\\
7	5\\
7	7\\
7	10\\
8	1\\
8	2\\
8	4\\
8	8\\
8	10\\
9	1\\
9	2\\
9	9\\
10	1\\
10	2\\
10	4\\
10	5\\
10	6\\
10	7\\
10	8\\
10	10\\
};
\end{axis}
\end{tikzpicture}%
       \caption{Original pattern.\label{fig:original_pattern}}
   \end{subfigure}
   \hfil
   \begin{subfigure}[t]{0.24\textwidth}
       \centering
\begin{tikzpicture}

\begin{axis}[%
width=0.8\textwidth,
height=0.8\textwidth,
at={(1.236in,0.481in)},
scale only axis,
xmin=0,
xmax=11,
xlabel style={font=\color{white!15!black}},
xlabel={nz = 52},
y dir=reverse,
ymin=0,
ymax=11,
axis background/.style={fill=white}
]
\addplot [color=red, only marks, mark size=2.3pt, mark={square*}, mark options={solid, red}, forget plot]
 table[row sep=crcr]{%
1	10\\
10	1\\
};
\addplot [color=blue, only marks, mark size=2.3pt, mark=*, mark options={solid, blue}, forget plot]
 table[row sep=crcr]{%
1	1\\
1	6\\
1	7\\
1	8\\
1	9\\
2	2\\
2	3\\
2	5\\
2	8\\
2	9\\
2	10\\
3	2\\
3	3\\
3	5\\
3	6\\
4	4\\
4	5\\
4	6\\
4	8\\
4	10\\
5	2\\
5	3\\
5	4\\
5	5\\
5	6\\
5	7\\
5	10\\
6	1\\
6	3\\
6	4\\
6	5\\
6	6\\
6	10\\
7	1\\
7	5\\
7	7\\
7	10\\
8	1\\
8	2\\
8	4\\
8	8\\
8	10\\
9	1\\
9	2\\
9	9\\
10	2\\
10	4\\
10	5\\
10	6\\
10	7\\
10	8\\
10	10\\
};
\end{axis}
\end{tikzpicture}%
       \caption{Modified pattern, with %
    \texttt{tol}$=10^{-3}$.\label{fig:removed_tol1}}
   \end{subfigure}
   \hfil
   \begin{subfigure}[t]{0.24\textwidth}
       \centering
\begin{tikzpicture}

\begin{axis}[%
width=0.8\textwidth,
height=0.8\textwidth,
at={(1.236in,0.481in)},
scale only axis,
xmin=0,
xmax=11,
xlabel style={font=\color{white!15!black}},
xlabel={nz = 50},
y dir=reverse,
ymin=0,
ymax=11,
axis background/.style={fill=white}
]
\addplot [color=green, only marks, mark size=2.3pt, mark={square*}, mark options={solid, green}, forget plot]
 table[row sep=crcr]{%
5	10\\
10	5\\
};
\addplot [color=blue, only marks, mark size=2.3pt, mark=*, mark options={solid, blue}, forget plot]
 table[row sep=crcr]{%
1	1\\
1	6\\
1	7\\
1	8\\
1	9\\
2	2\\
2	3\\
2	5\\
2	8\\
2	9\\
2	10\\
3	2\\
3	3\\
3	5\\
3	6\\
4	4\\
4	5\\
4	6\\
4	8\\
4	10\\
5	2\\
5	3\\
5	4\\
5	5\\
5	6\\
5	7\\
6	1\\
6	3\\
6	4\\
6	5\\
6	6\\
6	10\\
7	1\\
7	5\\
7	7\\
7	10\\
8	1\\
8	2\\
8	4\\
8	8\\
8	10\\
9	1\\
9	2\\
9	9\\
10	2\\
10	4\\
10	6\\
10	7\\
10	8\\
10	10\\
};
\end{axis}
\end{tikzpicture}%
       \caption{Modified pattern, with %
       \texttt{tol}$=10^{-6}$\label{fig:removed_tol2}.}
   \end{subfigure}
   \caption{Example of application of Algorithm~\ref{alg:adaptive_riemm} on a $10 \times 10$ reversible sparse matrix. We plot in \textcolor{red}{$\blacksquare$}  and \textcolor{green}{$\blacksquare$} the entries removed when \texttt{tol}$=10^{-3}$ and \texttt{tol}$=10^{-6}$, respectively.  }\label{fig:patter_removal}
\end{figure}

\section{Numerical examples}\label{sec:numerical-examples}

The implementation of the Riemannian optimizer relies on \texttt{Manopt}~\cite{manopt} (v.8.0), and it includes an implementation of the first order geometry for the manifolds $\mathcal{M}_{\boldsymbol{\pi}}$ in \eqref{eq:manifold_sparse_reversible} and $\mathcal{M}_{P,\boldsymbol{\pi}}$ in \eqref{eq:manifold_fixed_P}. The constrained version of the optimizer relies on the \texttt{fmincon} MATLAB function~\cite{MR1724768}. %
{The semidefinite programming approach discussed for comparison in Section~\ref{sec:semidefinite_approach} is implemented using the CVX library~\cite{cvx,gb08}.}
The code to replicate the experiments is available on the GitHub repository \url{https://www.github.com/Cirdans-Home/optimize-kemeny}.
Numerical experiments have been run on MATLAB R2022b, on a laptop with 16 GB of RAM and an Intel Core i7-10750H CPU (2.59 GHz).

\subsection{Comparing the Constrained and the Riemannian approaches}

Our comparison addresses three settings: (i) nearly reducible Markov chains~\cite{MR997458}, (ii) synthetic tests (including a scalability study as the state space grows), and (iii) reversible chains arising from power networks~\cite{MEDJROUBI201714}. The quality of the computed matrix $X$ is evaluated using 
$\|X\mathbf{1}-\mathbf{1}\|_{\infty}$, 
$\|\boldsymbol{\pi}^{\top}X-\boldsymbol{\pi}^{\top}\|_{\infty}$, 
and 
$\|D_{\boldsymbol{\pi}}X - X^{\top}D_{\boldsymbol{\pi}}\|_{\infty}$,
which measure, respectively, stochasticity, stationarity, and reversibility. We also report the Frobenius norm of the perturbation, running time (in seconds), and the change in Kemeny’s constant.

\begin{example}
Consider a reversible stochastic matrix $P$ associated with a nearly reducible Markov chain, obtained by constructing a $P$ as a stochastic matrix with random entries on the two blocks on the main diagonal and then adding extra small transition probabilities to the extra-diagonal blocks to make it irreducible. 
In this example, we consider a $50 \times 50$ test matrix $P$, with two diagonal blocks of size $25\times 25$ each, and two additional entries in positions $P_{1,50}$ and $P_{50,1}$. To make it  reversible, the Metropolis--Hasting modification described in Proposition~\ref{pro:metropolis-mod} is applied. %
We run both the constrained and the Riemannian optimizers on $P$. Given $\boldsymbol{\pi}$ the stationary probability vector of $P$, the solution we look for is a stochastic matrix $X$ with the same pattern of $P$ and such that $D_{\boldsymbol{\pi}} X = X^{\top} D_{\boldsymbol{\pi}}$.
Kemeny's constant associated with the original stochastic matrix $P$ is $2.2147 \times 10^{6}$. We test both the constrained and the Riemannian optimizers, comparing the quality of the solution. In Table~\ref{tab:Comparison_almost_red}, we report the different evaluation metrics. We observe that both approaches succeed in substantially reducing Kemeny's constant. Nevertheless, it is worth noticing that, for this example, the Riemannian optimizer, based on the Riemannian CG, is faster than the constrained one, without losing accuracy on the quality of the solution $X$. 

\begin{table}[htb]
\caption{Results for the nearly reducible matrix $P$.}
\label{tab:Comparison_almost_red}
\centering
\begin{tabular}{lcc}
\toprule
 & \textbf{Riemannian} & \textbf{Constrained} \\
\midrule
$\mathcal{K}(X)$ 
& $1.4746 \times 10^{2}$ 
& $1.0928\times 10^{2}$ \\

$\| X \mathbf{1} -\mathbf{1}\|_{\infty}$ 
& $6.6613\times 10^{-16}$ 
& $5.2637 \times 10^{-2}$ \\

$\| \boldsymbol{\pi}^{\top} X -\boldsymbol{\pi}^{\top}\|_{\infty}$ 
& $2.0816\times 10^{-17}$ 
& $1.8301\times10^{-3}$ \\

$\| D_{\boldsymbol{\pi}} X - X^{\top} D_{\boldsymbol{\pi}}\|_{\infty}$ 
& $4.9500\times 10^{-18}$ 
& $1.8559\times 10^{-3}$ \\

$\| X - P \|_F$ 
& $3.2895$ 
& $1.9479 \times 10^{-1}$ \\

Time (\si{\second}) 
& $2.6$ 
& $35.94$ \\

\bottomrule
\end{tabular}
\end{table}

\end{example}
\begin{example}
\label{ex:example_random}
    To compare the solvers, we consider a series of randomly generated problems. %
After choosing a dimension $n$, we construct a test reversible stochastic matrix with a specified pattern $\mathcal{P}$ and choosing $\mathcal{S} \equiv \mathcal{P}$ for the modification pattern. %
Our implementation of a first order geometry for the manifold \eqref{eq:manifold_sparse_reversible} is contained in the function  \lstinline[style=MATLAB-editor,breaklines=true]{N = multinomialsparsesymmetricfixedfactory(pv,S)}.

    In particular, for each choice of the size $n\in \left\lbrace 10,20,30,40,50,60 \right\rbrace$, we generate $25$ test cases, ordered with respect to the dimension $n$, and compare three solvers: the constrained one, the Riemannian approach using as optimizer the Riemannian CG~\cite{MR3325229}, and the Riemannian approach employing BB line search rule~\cite{Iannazzo}. In the subsequent figures, on the $x$-axis the test cases are numbered, $25$ for each dimension, for the dimensions $n\in \left\lbrace 10,20,30,40,50,60 \right\rbrace$.

    \begin{figure}[htb]
        \centering
        \begin{tikzpicture}
            \begin{axis}[
            legend pos= south east,
            width=\linewidth,
            height=45mm,
            xmin=1, xmax= 150,
            xtick={25,50,75,100,125,150},
              xlabel={Test cases},
            xlabel style={yshift=0.5em},
            legend style={draw=none,fill=none,legend cell align=left, align=left, font=\footnotesize}
            ]
                \addplot[only marks, cyan, mark=x, mark size = 2.5pt] table[y index = 1] {kem_comparison_s.dat};
                \addplot[only marks, blue, mark=triangle, mark size = 2.5pt] table[y index = 2] {kem_comparison_s.dat};
               \addplot[only marks, color=mycolor2, mark options={solid, mycolor2}, mark=o, mark size = 2.5pt] table[y index = 3] {kem_comparison_s.dat};
                \addplot[color=mycolor3, only marks, mark=square, mark options={solid, mycolor3}, mark size = 2.5pt] table[y index = 4] {kem_comparison_s.dat};
                \legend{$\mathcal{K}(P)$, $\mathcal{K}(X)$ with Riem-CG, $\mathcal{K}(X)$ with Riem-BB,  $\mathcal{K}(X)$ with Constr.}
            \end{axis}
        \end{tikzpicture}
        \caption{Kemeny's constant of $X$ and of the original matrix $P$. %
        }
        \label{fig:comparison_wrt_kem}
    \end{figure}
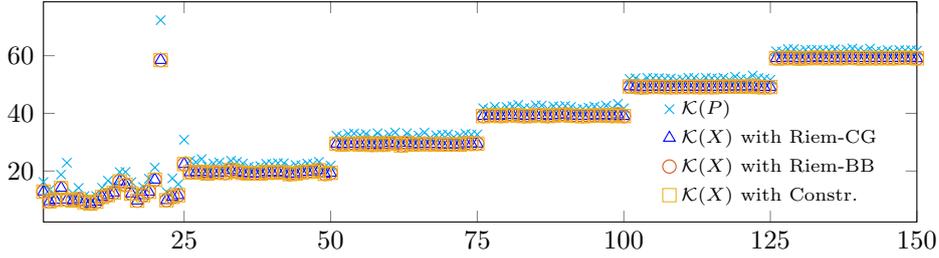
    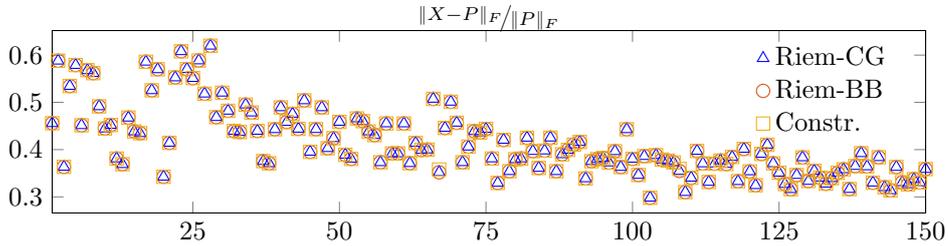
\begin{figure}[htb]
        \centering
        \begin{tikzpicture}
            \begin{axis}[
            legend pos= north east,
            width=\linewidth,
            height=40mm,
            xmin=1, xmax= 150,
            xtick={25,50,75,100,125,150},
              xlabel={Test cases},
            xlabel style={yshift=0.5em},
            legend style={draw=none,fill=none,legend cell align=left, align=left, font=\normalsize},
            title={$\nicefrac{\| X- P \|_F}{\|P\|_F}$},
            title style={yshift=-1em},
            ]
                \addplot[only marks, blue, mark=triangle, mark size = 2.5pt] table[y index = 1] {dist_comparison.dat};
               \addplot[only marks, color=mycolor2, mark options={solid, mycolor2}, mark=o, mark size = 2.5pt] table[y index = 2] {dist_comparison.dat};
                \addplot[color=mycolor3, only marks, mark=square, mark options={solid, mycolor3}, mark size = 2.5pt] table[y index = 3] {dist_comparison.dat};
            \legend{Riem-CG,  Riem-BB, Constr.}
            \end{axis}
        \end{tikzpicture}
        \caption{Detailed comparison of the relative distance for each test problem.} %
        \label{fig:comparison_wrt_dist}
    \end{figure}
    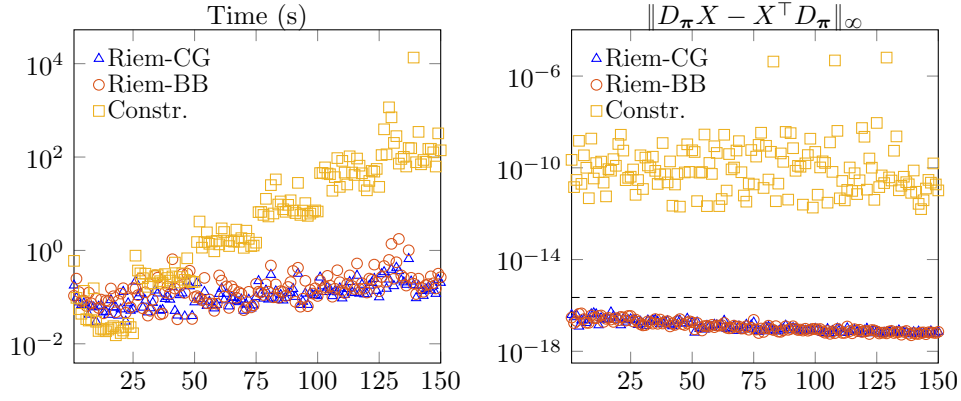
\begin{figure}[htb]
        \centering    \begin{tikzpicture}
            \begin{semilogyaxis}[
           legend style={nodes={scale=0.6, transform shape}},
           legend style={font=\Large},
            legend pos= north west,
            width=0.49\linewidth,
            height=60mm,
            xmin=1, xmax= 150,
            xtick={25,50,75,100,125,150},
              xlabel={Test cases},
            xlabel style={yshift=0.5em},
          legend style={draw=none,fill=none,legend cell align=left, align=left},
          title={Time (s)},
          title style={yshift=-1em},
            ]
                \addplot[only marks, blue, mark=triangle] table[y index = 1] {times_comparison.dat};
               \addplot[only marks, color=mycolor2, mark options={solid, mycolor2}, mark=o] table[y index = 2] {times_comparison.dat};
                \addplot[color=mycolor3, only marks, mark=square, mark options={solid, mycolor3}] table[y index = 3] {times_comparison.dat};
                \legend{Riem-CG,  Riem-BB, Constr.}
            \end{semilogyaxis}
        \end{tikzpicture}~
        \begin{tikzpicture}
            \begin{semilogyaxis}[
             legend style={nodes={scale=0.6, transform shape}},
           legend style={font=\Large},
            legend pos= north west,
            width=0.49\linewidth,
            height=60mm,
            xmin=1, xmax= 150,
            xtick={25,50,75,100,125,150},
              xlabel={Test cases},
            xlabel style={yshift=0.5em},
              legend style={draw=none,fill=none,legend cell align=left, align=left},
              title={$\|D_{\boldsymbol{\pi}} X - X^{\top} D_{\boldsymbol{\pi}}\|_{\infty}$},
              title style={yshift=-1em},
            ]
                \addplot[only marks, blue, mark=triangle] table[y index = 1] {rev_comparison.dat};
               \addplot[color=mycolor2, mark options={solid, mycolor2}, only marks, mark=o] table[y index = 2] {rev_comparison.dat};
                \addplot[color=mycolor3, only marks, mark=square, mark options={solid, mycolor3}] table[y index = 3] {rev_comparison.dat};
                  \addplot [color=black, dashed, forget plot]
                  table[row sep=crcr]{%
                1	2.22044604925031e-16\\
                150	2.22044604925031e-16\\
                };BB
                \legend{Riem-CG,  Riem-BB, Constr.}
            \end{semilogyaxis}
        \end{tikzpicture}
        \caption{Comparison with respect to the computational time, and with respect to the metric $\| D_{\boldsymbol{\pi}} X - X^{\top} D_{\boldsymbol{\pi}}\|_{\infty}$. %
        }
        \label{fig:comparison_wrt_rev_time}
    \end{figure}
      \begin{figure}[htb]
        \centering
        \begin{tikzpicture}
            \begin{semilogyaxis}[
             legend style={nodes={scale=0.6, transform shape}},
           legend style={font=\Large},
            legend pos= north west,
            width=0.49\linewidth,
            height=60mm,
            xmin=1, xmax= 150,
            xtick={25,50,75,100,125,150},
            xlabel={Test cases},
            xlabel style={yshift=0.5em},
      legend style={draw=none,fill=none,legend cell align=left, align=left},
      title ={$\|\boldsymbol{\pi}^{\top}X- \boldsymbol{\pi}^{\top}\|_{\infty}$},
      title style={yshift=-1em},
            ]
                \addplot[only marks, blue, mark=triangle] table[y index = 1] {stat_comparison.dat};
               \addplot[only marks, color=mycolor2, mark options={solid, mycolor2}, mark=o] table[y index = 2] {stat_comparison.dat};
                \addplot[color=mycolor3, only marks, mark=square, mark options={solid, mycolor3}] table[y index = 3] {stat_comparison.dat};
                \addplot [color=black, dashed, forget plot]
                  table[row sep=crcr]{%
                1	2.22044604925031e-16\\
                150	2.22044604925031e-16\\
                };
                                \legend{Riem-CG,  Riem-BB, Constr.}
            \end{semilogyaxis}
        \end{tikzpicture}~
         \begin{tikzpicture}
            \begin{semilogyaxis}[
             legend style={nodes={scale=0.6, transform shape}},
           legend style={font=\Large},
            legend pos= north west,
            width=0.49\linewidth,
            height=60mm,
            xmin=1, xmax= 150,
            xtick={25,50,75,100,125,150},
              xlabel={Test cases},
            xlabel style={yshift=0.5em},
                  legend style={draw=none,fill=none,legend cell align=left, align=left},
                  title ={$\|X\mathbf{1} - \mathbf{1}\|_{\infty}$},
                  title style={yshift=-1em},
            ]
                \addplot[only marks, blue, mark=triangle] table[y index = 1] {stoc_comparison.dat};
               \addplot[only marks, color=mycolor2, mark options={solid, mycolor2}, mark=o] table[y index = 2] {stoc_comparison.dat};
                \addplot[color=mycolor3, only marks, mark=square, mark options={solid, mycolor3}] table[y index = 3] {stoc_comparison.dat};
                  \addplot [color=black, dashed, forget plot]
                  table[row sep=crcr]{%
                1	2.22044604925031e-16\\
                150	2.22044604925031e-16\\
                };
            \legend{Riem-CG,  Riem-BB, Constr.}
            \end{semilogyaxis}
        \end{tikzpicture}
        \caption{Comparison with respect to the metrics $\| \boldsymbol{\pi}^{\top} X -\boldsymbol{\pi}^{\top} \|_{\infty}$ and $\| X\mathbf{1} -\mathbf{1}\|_{\infty}$. %
        }
        \label{fig:comparison_wrt_stat_stoc}
    \end{figure}
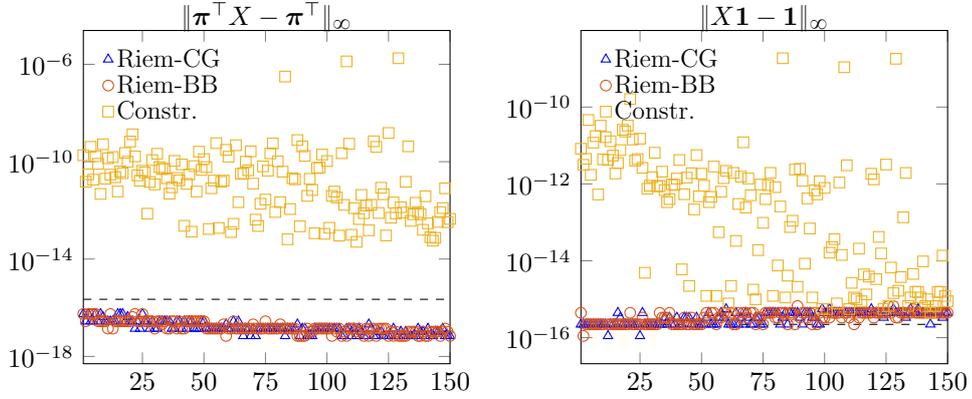

From Figure \ref{fig:comparison_wrt_kem} and \ref{fig:comparison_wrt_dist}, we observe that the solution $X$ computed via the constrained approach and the Riemannian one---for both choices of the inner Riemannian optimizer---produce the same reduction of the Kemeny's constant and relative distance $\nicefrac{\| X- P \|_{F}}{\|P \|_F}$. 

The main difference between the constrained version of the approach and the Riemannian one consists in the quality of the produced solution, as it can be seen from Figures \ref{fig:comparison_wrt_rev_time}(Right) and \ref{fig:comparison_wrt_stat_stoc}. Our implementation of the constrained optimizer employs \texttt{fmincon} function in MATLAB, which appears to suffer in this set of test problems. Moreover, a second difference may be found in the computational time employed by the three methods, as displayed in Figure \ref{fig:comparison_wrt_rev_time}(Left). Indeed, while the constrained approach seems faster for small sized stochastic matrices, the Riemannian versions are faster for large sized test cases.
\end{example}
\begin{example}
In an analogous way to the procedure in Example \ref{ex:example_random}, we compare the three solvers on a set of randomly generated reversible stochastic matrices, where we choose a pattern $\mathcal{S}$ different from the pattern $\mathcal{P}$ of the original reversible matrix $P$. The comparison has been carried out with respect to the minimization of Kemeny's constant (Fig.~\ref{fig:fixed_kemeny}), the relative distance (Fig.~\ref{fig:fixed_wrt_dist}), to the elapsed time in Fig.~\ref{fig:fixed_wrt_rev_time} (Left) and the metrics in Fig.~\ref{fig:fixed_wrt_rev_time} (Right) and Fig.~\ref{fig:fixed_wrt_stat_stoc}. The test set can be generated using the implementation of the manifold $\mathcal{M}_{P,\boldsymbol{\pi}}$, stored in 
\lstinline[style=MATLAB-editor]{N=multinomialsymmetricfixedentriesfactory(pv,S,calP)}, available on Github.

%
For each choice of the size in $n\in \left\lbrace 10,20,30,40,50,60\right\rbrace$, we generate $25$ test matrices and order them with respect to the dimension. In each figure, on the $x$-axis the test cases are numbered, following this order.

Similarly to the results in Example~\ref{ex:example_random}, the Riemannian solvers have a comparable behaviour, while the constrained optimizer seems to suffer when the dimension increases, leading to a loss in the accuracy of the produced solutions and a higher computational time. However, for lower dimensions, the constrained procedure provides a valuable alternative to the Riemannian optimizer. Observe that both Riemannian solvers do not suffer when dealing with matrices of higher dimensions and, more importantly, the quality of the produced solution does not depend on the size of the matrices involved.

       \begin{figure}
        \centering
        \begin{tikzpicture}
            \begin{axis}[
            legend pos= south east,
            width=\linewidth,
            height=45mm,
            xmin=1, xmax= 150,
            xlabel={Test cases},
            xlabel style={yshift=0.5em},
            xtick={25,50,75,100,125,150},
            legend style={draw=none,fill=none,legend cell align=left, align=left, font=\footnotesize}
            ]
                \addplot[only marks, cyan, mark=x, mark size = 2.5pt] table[y index = 1] {kem_fixed_s.dat};
                \addplot[only marks, blue, mark=triangle, mark size = 2.5pt] table[y index = 2] {kem_fixed_s.dat};
               \addplot[only marks, color=mycolor2, mark options={solid, mycolor2}, mark=o, mark size = 2.5pt] table[y index = 3] {kem_fixed_s.dat};
               \addplot[color=mycolor3, only marks, mark=square, mark options={solid, mycolor3}, mark size = 2.5pt] table[y index = 4] {kem_fixed_s.dat};
                \legend{$\mathcal{K}(P)$, $\mathcal{K}(X)$ with Riem-CG, $\mathcal{K}(X)$ with Riem-BB,  $\mathcal{K}(X)$ with Constr.}
            \end{axis}
        \end{tikzpicture}
        \caption{Kemeny's constant of $X$ and the original matrix $P$. %
        }
        \label{fig:fixed_kemeny}
    \end{figure}
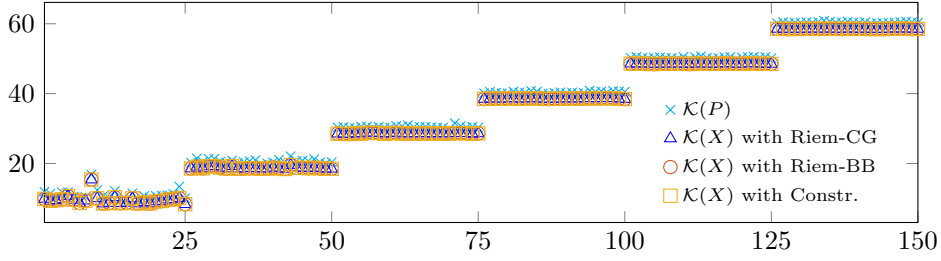

    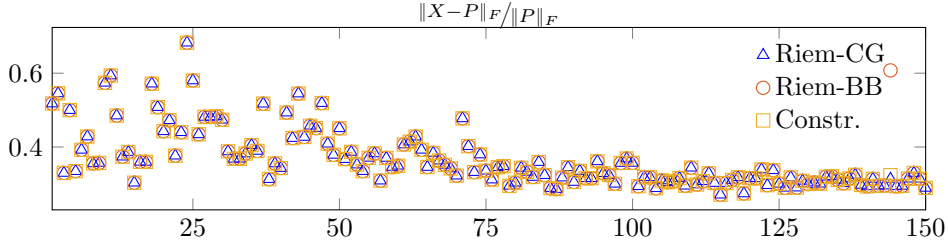
\begin{figure}
        \centering
        \begin{tikzpicture}
            \begin{axis}[
            legend pos= north east,
            width=\linewidth,
            height=40mm,
            xmin=1, xmax= 150,
            xtick={25,50,75,100,125,150},
            xlabel={Test cases},
            xlabel style={yshift=0.5em},
            title={$\nicefrac{\| X- P \|_F}{\|P\|_F}$},
            title style={yshift=-1em},
            legend style={draw=none,fill=none,legend cell align=left, align=left}
            ]
                \addplot[only marks, blue, mark=triangle, mark size = 2.5pt] table[y index = 1] {dist_fixed.dat};
               \addplot[only marks, color=mycolor2, mark options={solid, mycolor2}, mark=o, mark size = 2.5pt] table[y index = 2] {dist_fixed.dat};
                \addplot[color=mycolor3, only marks, mark=square, mark options={solid, mycolor3}, mark size = 2.5pt] table[y index = 3] {dist_fixed.dat};
                \legend{Riem-CG, Riem-BB, Constr.}
            \end{axis}
        \end{tikzpicture}
       \caption{Detailed comparison of the relative distance, for each test problem. %
       }
       \label{fig:fixed_wrt_dist}
    \end{figure}

    \begin{figure}
        \centering    \begin{tikzpicture}
            \begin{semilogyaxis}[
            legend style={nodes={scale=0.6, transform shape}},
           legend style={font=\Large},
            legend pos= north west,
            width=0.49\linewidth,
            height=45mm,
            xmin=1, xmax= 150,
            xtick={25,50,75,100,125,150},
            xlabel={Test cases},
            xlabel style={yshift=0.5em},
          legend style={draw=none,fill=none,legend cell align=left, align=left},
          title ={Time (s)},
          title style={yshift=-1em},
            ]
                \addplot[only marks, blue, mark=triangle] table[y index = 1] {time_fixed.dat};
               \addplot[only marks, color=mycolor2, mark options={solid, mycolor2}, mark=o] table[y index = 2] {time_fixed.dat};
               \addplot[color=mycolor3, only marks, mark=square, mark options={solid, mycolor3}] table[y index = 3] {time_fixed.dat};
               \legend{Riem-CG, Riem-BB, Constr.}
            \end{semilogyaxis}
        \end{tikzpicture}~
        \begin{tikzpicture}
            \begin{semilogyaxis}[
             legend style={nodes={scale=0.6, transform shape}},
           legend style={font=\Large},
            legend pos= north west,
            width=0.49\linewidth,
            height=45mm,
            xmin=1, xmax= 150,
            xtick={25,50,75,100,125,150},
              xlabel={Test cases},
            xlabel style={yshift=0.5em},
              legend style={draw=none,fill=none,legend cell align=left, align=left},
              title ={$\|D_{\boldsymbol{\pi}} X - X^{\top} D_{\boldsymbol{\pi}}\|_{\infty}$},
              title style={yshift=-1em},
            ]
                \addplot[only marks, blue, mark=triangle] table[y index = 1] {rev_fixed.dat};
               \addplot[color=mycolor2, mark options={solid, mycolor2}, only marks, mark=o] table[y index = 2] {rev_fixed.dat};
               \addplot[color=mycolor3, only marks, mark=square, mark options={solid, mycolor3}] table[y index = 3] {rev_fixed.dat};
                  \addplot [color=black, dashed, forget plot]
                  table[row sep=crcr]{%
                1	2.22044604925031e-16\\
                150	2.22044604925031e-16\\
                };
                               \legend{Riem-CG, Riem-BB, Constr.}
            \end{semilogyaxis}
        \end{tikzpicture}
        \caption{Comparison with respect to the computational time, and with respect to the metric $\| D_{\boldsymbol{\pi}} X - X^{\top} D_{\boldsymbol{\pi}}\|_{\infty}$. %
        }
        \label{fig:fixed_wrt_rev_time}
    \end{figure}
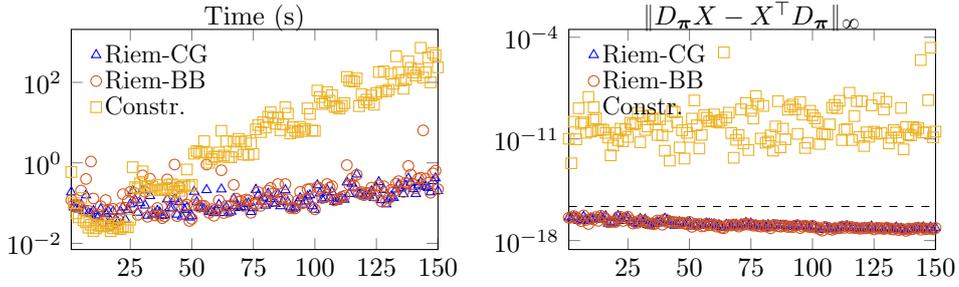

      \begin{figure}
        \centering
        \begin{tikzpicture}
            \begin{semilogyaxis}[
             legend style={nodes={scale=0.6, transform shape}},
           legend style={font=\Large},
            legend pos= north west,
            width=0.49\linewidth,
            height=45mm,
            xmin=1, xmax= 150,
            xtick={25,50,75,100,125,150},
            xlabel={Test cases},
            xlabel style={yshift=0.5em},
      legend style={draw=none,fill=none,legend cell align=left, align=left},
      title ={$\|\boldsymbol{\pi}^{\top}X- \boldsymbol{\pi}^{\top}\|_{\infty}$},
      title style={yshift=-1em},
            ]
                \addplot[only marks, blue, mark=triangle] table[y index = 1] {stat_fixed.dat};
               \addplot[only marks, color=mycolor2, mark options={solid, mycolor2}, mark=o] table[y index = 2] {stat_fixed.dat};
                \addplot[color=mycolor3, only marks, mark=square, mark options={solid, mycolor3}] table[y index = 3] {stat_fixed.dat};
                \addplot [color=black, dashed, forget plot]
                  table[row sep=crcr]{%
                1	2.22044604925031e-16\\
                150	2.22044604925031e-16\\
                };
           \legend{Riem-CG, Riem-BB, Constr.}
            \end{semilogyaxis}
        \end{tikzpicture}~
         \begin{tikzpicture}
            \begin{semilogyaxis}[
            legend style={nodes={scale=0.6, transform shape}},
           legend style={font=\Large},
            legend pos= north west,
            width=0.49\linewidth,
            height=45mm,
            xmin=1, xmax= 150,
            xtick={25,50,75,100,125,150},
            xlabel={Test cases},
            xlabel style={yshift=0.5em},
                  legend style={draw=none,fill=none,legend cell align=left, align=left},
            title={$\|X\mathbf{1} - \mathbf{1}\|_{\infty}$},
            title style={yshift=-1em},
            ]
                \addplot[only marks, blue, mark=triangle] table[y index = 1] {stoc_fixed.dat};
               \addplot[only marks, color=mycolor2, mark options={solid, mycolor2}, mark=o] table[y index = 2] {stoc_fixed.dat};
               \addplot[color=mycolor3, only marks, mark=square, mark options={solid, mycolor3}] table[y index = 3] {stoc_fixed.dat};
                  \addplot [color=black, dashed, forget plot]
                  table[row sep=crcr]{%
                1	2.22044604925031e-16\\
                150	2.22044604925031e-16\\
                };
                \legend{Riem-CG, Riem-BB, Constr.}
            \end{semilogyaxis}
        \end{tikzpicture}
        \caption{Comparison with respect to the metrics $\| \boldsymbol{\pi}^{\top} X -\boldsymbol{\pi}^{\top} \|_{\infty}$ and $\| X\mathbf{1} -\mathbf{1}\|_{\infty}$. %
        }
        \label{fig:fixed_wrt_stat_stoc}
    \end{figure}
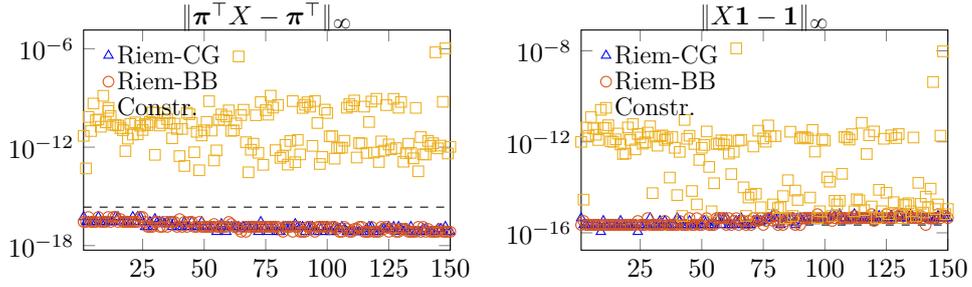

\end{example}

\begin{example}
We consider a set of the power grid networks of several countries taken from the dataset \texttt{Power grid}\footnote{The dataset is available at \url{https://github.com/ComplexNetTSP/Power_grids/tree/v1.0.0}}, which provides the adjacency matrices of the networks; see~\cite{MEDJROUBI201714} for a description of the extraction procedure for these data.  We select $5$ power networks, reducing the problem to their largest connected component when needed. For each reduced network with adjecency matrix $A$, we may then construct a random walk as in Definition \ref{def:random-walk}, scaling the entries as:
$P = D_{\mathbf{d}}^{-1}A$, with $\mathbf{d} = A\mathbf{1}.$
The matrix $P$ is irreducible and, furthermore, the stationary vector for $P$ is $\boldsymbol{\pi}=\frac{\mathbf{d}}{\sum_i \mathbf{d}_i}$. We observe that the transition matrix $P$ is reversible with respect to $\boldsymbol{\pi}$.

We may then apply the constrained optimizer and the Riemannian optimizer based on the Riemannian-CG on the reversible matrix $P$, choosing as pattern the one of $P$.
%
%
%
In Table~\ref{tab:comparison power grid}, we provide the results for the final matrix $X$, obtained via both the constrained and the Riemannian optimizer, comparing its Kemeny's constant, the reversibility error $\| D_{\boldsymbol{\pi}} X - X^{\top} D_{\boldsymbol{\pi}}\|_{\infty}$ and the elapsed time in seconds. Moreover, in Table~\ref{tab:summary of power grids} we summarize the features of the considered test cases, reporting for each the matrix size $n$, the number of nonzero elements, the density (computed as the percentage of nonzero entries out of $n^2$), along with its original Kemeny's constant $\mathcal{K}(P)$ and the bound in~\eqref{eq:kirkland-bound}.
All quantities refer to the largest connected component, which is the one used in the experiments---we remark that the power network associated with the Switzerland is irreducible. While the reduction of Kemeny's constant is comparable, the Riemannian optimizer satisfies the reversibility constraints with better accuracy.

\begin{table}[htbp]
\caption{Structural properties and original Kemeny's constants for the networks in \texttt{Power grid}.}
\label{tab:summary of power grids}
\centering
\small

\begin{tabular}{l c c c c c}
\toprule
Network
& Size
& $\operatorname{nnz}(P)$
& Density
& $\mathcal{K}(P)$
& Bound~\eqref{eq:kirkland-bound}
\\

\midrule

Austria
& $147$
& $336$
& $1.55\%$
& $1.1896\times10^{3}$
& $9.3027\times10^{1}$
\\

Belgium
& $90$
& $218$
& $2.69\%$
& $5.2136\times10^{2}$
& $5.7234\times10^{1}$
\\

Denmark
& $63$
& $136$
& $3.42\%$
& $6.6786\times10^{2}$
& $3.8342\times10^{1}$
\\

Netherlands
& $84$
& $190$
& $2.69\%$
& $5.2300\times10^{2}$
& $5.2191\times10^{1}$
\\

Switzerland
& $310$
& $736$
& $0.76\%$
& $1.9366\times10^{4}$
& $1.9748\times10^{2}$
\\

\bottomrule
\end{tabular}
\end{table}

\begin{table}[htbp]
\caption{Performance comparison for Markov chains associated with networks in \texttt{Power grid}.}
\label{tab:comparison power grid}
\centering
\small
\setlength{\tabcolsep}{4pt}

\begin{tabular}{l ccc ccc}
\toprule
&
\multicolumn{3}{c}{Constrained}
&
\multicolumn{3}{c}{Riemannian}
\\

\cmidrule(lr){2-4}
\cmidrule(lr){5-7}

Network
& $\mathcal{K}(X)$
& {\footnotesize $\left\| D_{\boldsymbol{\pi}} X - X^{\top} D_{\boldsymbol{\pi}} \right\|_{\infty}$} %
& Time(\si{\second})
& $\mathcal{K}(X)$
& {\footnotesize $\left\| D_{\boldsymbol{\pi}} X - X^{\top} D_{\boldsymbol{\pi}} \right\|_{\infty}$} %
& Time(\si{\second})
\\

\midrule

Austria
& $1.1550\times10^{3}$
& $1.6067\times10^{-5}$
& $161.80$
& $1.1535\times10^{3}$
& $5.6379\times10^{-18}$
& $17.91$
\\

Belgium
& $4.9493\times10^{2}$
& $1.3619\times10^{-6}$
& $38.62$
& $4.9482\times10^{2}$
& $3.4694\times10^{-18}$
& $6.12$
\\

Denmark
& $6.5814\times10^{2}$
& $1.6731\times10^{-5}$
& $9.36$
& $6.4610\times10^{2}$
& $9.5410\times10^{-18}$
& $4.59$
\\

Netherlands
& $5.0549\times10^{2}$
& $2.0134\times10^{-9}$
& $21.19$
& $5.0549\times10^{2}$
& $9.5410\times10^{-18}$
& $5.50$
\\

Switzerland
& $1.3900\times10^{4}$
& $1.9053\times10^{-7}$
& $40.16$
& $2.6217\times10^{3}$
& $2.3852\times10^{-18}$
& $102.98$
\\

\bottomrule
\end{tabular}
\end{table}

In Fig.~\ref{fig:switzerland} (left), we plot the original power graph of Switzerland. The application of one of the two optimizers produces a new Markov chain, with the same stationary probability vector $\boldsymbol{\pi}$. In this setting, it is possible to track the modifications to the original graph, constructing the new adjacency matrix as $\tilde{A} = D_{\mathbf{d}} X$. In Fig.~\ref{fig:switzerland} (right), we plot the final graph obtained via the Riemannian optimizer, highlighting the edges whose weight increases (blue) and the ones whose weight decreases after the minimization procedure (green).

\begin{figure}[htbp]
    \centering

    \begin{minipage}{0.54\textwidth}
        \includegraphics[width=\linewidth,height=0.6\textheight,keepaspectratio]{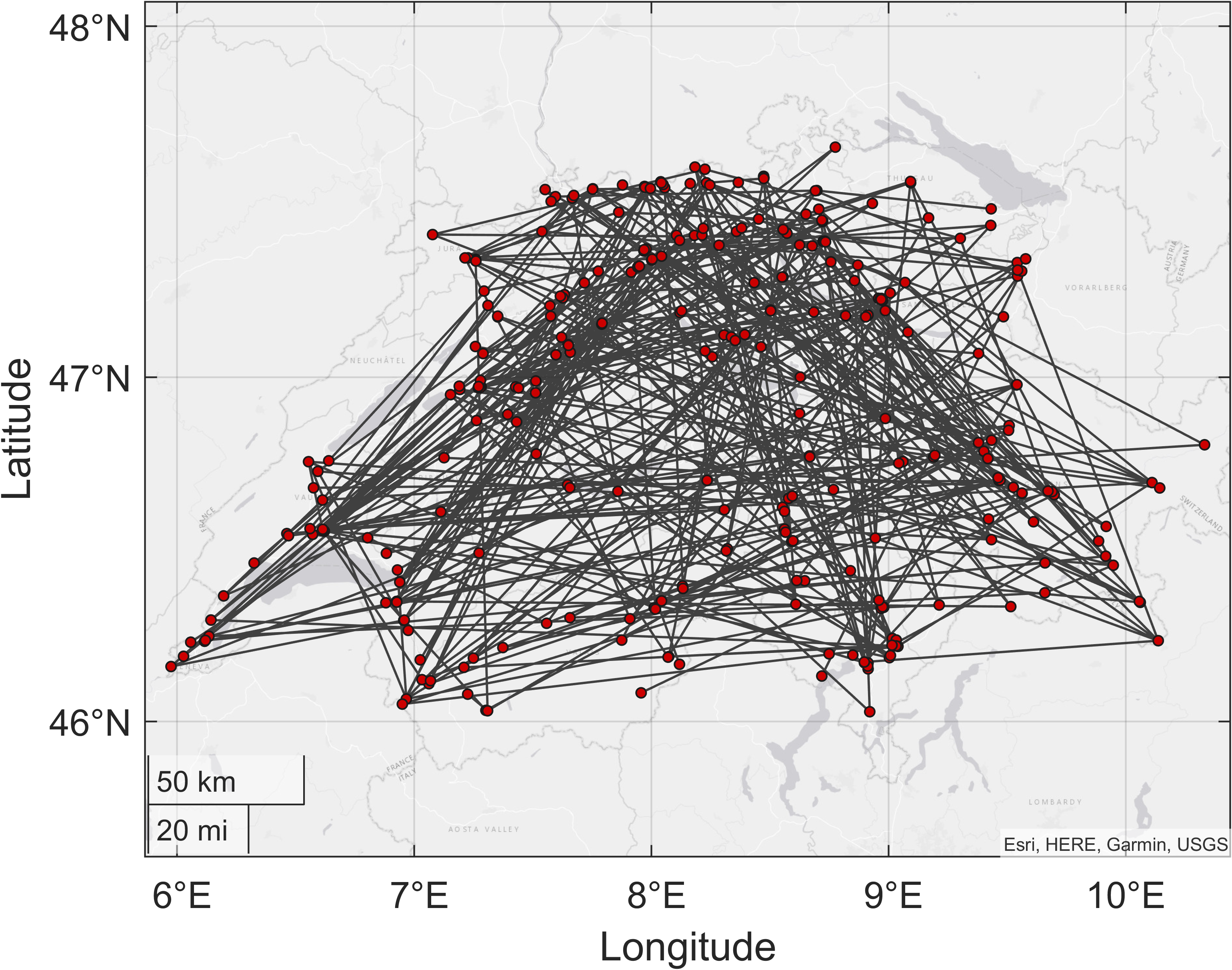}
    \end{minipage}%
    \hfill
    \begin{minipage}[t]{0.45\textwidth}
        \includegraphics[width=\linewidth,height=0.15\textheight,keepaspectratio]{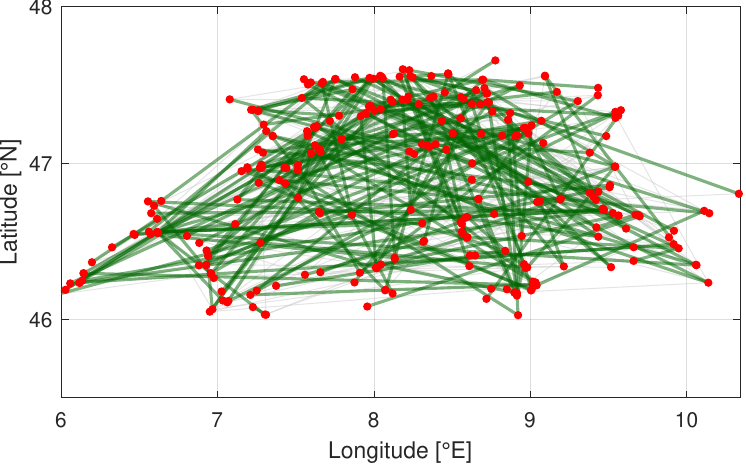}\\[0pt]
        \includegraphics[width=\linewidth,height=0.15\textheight,keepaspectratio]{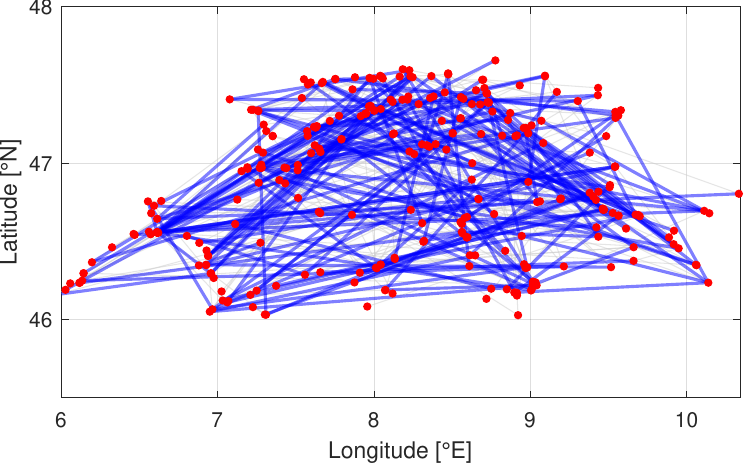}
    \end{minipage}

    \caption{Power network of Switzerland (on the left). Edges whose weights increased after optimizing Kemeny's constant are shown in blue, while those whose weights decreased are shown in green (on the right).}
    \label{fig:switzerland}
\end{figure}

\end{example}

\subsection{{Comparison with the semidefinite programming approach}}\label{sec:semidefinite_approach}

{Minimization of Kemeny’s constant has also been investigated in the context of stochastic surveillance strategies based on Markov chains~\cite{Bullo,BulloReview}. The perspective adopted in that setting differs substantially from ours. In particular, the problem considered in~\cite{Bullo,BulloReview} consists of determining, for a prescribed sparsity pattern of the transition matrix and a given stationary distribution, the realization of the transition probabilities that minimizes the Kemeny constant.} 

{By contrast, in our framework the objective is to determine a low-norm, fixed-pattern perturbation of an existing Markov chain that achieves a reduction of the Kemeny constant while preserving stochasticity and structural constraints.}

{The realization problem considered in~\cite{Bullo,BulloReview} can be reformulated, in our notation, as the following semidefinite programming (SDP) problem:
\begin{equation}
\begin{aligned}
\min_{X,P} \quad & \operatorname{tr}(X) \\
\text{subject to} \quad &
\begin{bmatrix}
I - D_{\hat{\boldsymbol{\pi}}} P D_{\hat{\boldsymbol{\pi}}}^{-1} + \hat{\boldsymbol{\pi}} \hat{\boldsymbol{\pi}}^{\top} & I \\
I & X
\end{bmatrix}
\succeq 0, \\
& P \mathbf{1} = \mathbf{1}, \\
& \pi_i P_{ij} = \pi_j P_{ji},
\qquad \{i,j\}\in \mathcal{S}, \\
& 0 \le P_{ij} \le 1,
\qquad \{i,j\}\in \mathcal{S}, \\
& P_{ij} = 0,
\qquad \{i,j\}\notin \mathcal{S},
\end{aligned}
\label{eq:kemeny_sdp}
\end{equation}
where the pattern $\mathcal{S}$ encodes the admissible transitions of the Markov chain, and we have used the notation $M \succeq 0$ to state that $M$ is positive semidefinite. The linear matrix inequality in~\eqref{eq:kemeny_sdp} guarantees reversibility and allows the minimization of the Kemeny constant through the trace objective $\operatorname{tr}(X)$. The resulting optimization problem is convex and can therefore be efficiently solved using standard SDP solvers. To showcase the difference with our implementation, the problem~\eqref{eq:kemeny_sdp} has been formulated and solved in MATLAB using the CVX package~\cite{cvx,gb08}.}

\begin{example}
{For the numerical experiments, we considered the same construction used for the synthetic test cases in Example~\ref{ex:example_random} reversible Markov chains of dimension $n=30$. A random stationary distribution was generated by sampling positive entries and normalizing them so that they define a probability vector. Starting from this distribution, we constructed random sparse transition matrices compatible with a prescribed graph topology. The sparsity pattern was obtained from a randomly generated symmetric adjacency matrix with self-loops, ensuring that the resulting Markov chain remained connected and admissible. A symmetric matrix having the prescribed stationary distribution as invariant eigenvector was then generated using the manifold construction implemented in the routine \texttt{multinomialsparsesymmetricfixedfactory}. This matrix was converted into a reversible stochastic matrix satisfying detailed balance with respect to the chosen stationary distribution. Finally, the Kemeny constant associated with the generated Markov chain was computed and used as the baseline reference value in the subsequent optimization experiments.}

{Starting from an initial reversible Markov chain with Kemeny constant equal to \(33.485261\), the proposed Riemannian optimization approach produced a perturbed chain with Kemeny constant \(30.366179\) in \(1.192\times 10^{-1}\,\si{\second}\). The obtained perturbation remained relatively close to the original chain, with relative distance \(4.179491\times 10^{-1}\). Moreover, the constraints defining reversibility, stochasticity, and stationarity were all satisfied up to machine precision, yielding residuals \(1.06\times 10^{-17}\), \(2.22\times 10^{-16}\), and \(1.39\times 10^{-17}\), respectively.}

{Next, the standard Euclidean constrained optimization framework produced a local minimum with a Kemeny constant of \(30.366190\), which is nearly identical to the value achieved by the Riemannian approach with an execution time \qty{1.2672}{\second}. However, the constraints were satisfied to a lower degree of precision compared to the manifold-based framework, yielding residuals of \(1.17\times 10^{-11}\), \(3.00\times 10^{-14}\), and \(2.18\times 10^{-12}\) for the reversibility, stochasticity, and stationarity conditions, respectively.}

{Finally, the SDP-based formulation achieved a slightly lower Kemeny constant, equal to \(30.28066\), but required a significantly larger modification of the original chain, with relative distance \(1.086508\). The computational cost was also higher, with execution time \(2.3623\,\si{\second}\), more than one order of magnitude larger than that of the Riemannian method. The resulting matrix still satisfied the stochasticity, stationarity, and reversibility constraints to high numerical accuracy, with residuals \(4.51\times 10^{-12}\), \(2.31\times 10^{-13}\), and \(1.15\times 10^{-17}\), respectively.}

{Overall, the experiments indicate that the SDP approach is slightly more effective in directly minimizing the Kemeny constant, whereas the proposed Riemannian framework achieves a comparable reduction while preserving significantly greater proximity to the original Markov chain.}
\end{example}

\section{Conclusions}\label{sec:conclusions}

In this work, we investigated the problem of perturbing the entries of the transition matrix of a homogeneous, discrete-time Markov chain in order to reduce its Kemeny's constant, thus improving its circulation properties. While the general formulation of this problem is inherently non-convex and does not admit systematic solution methods, we identified a tractable special case by imposing reversibility of the underlying Markov chain. This additional structure renders the minimization of the Kemeny's constant a convex problem, making it amenable to the use of optimization techniques.

We proposed two complementary approaches: one based on constrained optimization via an IPM, and the other on Riemannian optimization. Both were analyzed theoretically and translated into algorithms exploiting the problem's structure. Numerical experiments show that while the IPM handles sparse cases and naturally zeroed pattern elements effectively, it may struggle to satisfy machine-precision constraints, particularly reversibility and preservation of the stationary distribution. In contrast, the Riemannian approach is more robust in maintaining these structural properties.

Several directions remain for future work, including connections to related optimization problems in spectral graph theory and stochastic processes, as well as the development of a second-order geometry on the proposed manifolds to enable faster methods.

\backmatter

\bmhead{Acknowledgements}
We thank the anonymous referees for their valuable comments and suggestions, which helped improve the quality of this paper.

\section*{Declarations}

{\small\bmhead{Funding} All the authors are member of the INdAM GNCS and acknowledge the MUR Excellence Department Project awarded to the Department of Mathematics, University of Pisa, CUP I57G22000700001. This work was partially supported by the Italian Ministry of University and Research (MUR) through the PRIN 2022 ``Low-rank Structures and Numerical Methods in Matrix and Tensor Computations and their Application'' code: 20227PCCKZ MUR D.D. financing decree n. 104 of February 2nd, 2022 (CUP I53D23002280006) and through the PRIN 2022 ``MOLE: Manifold constrained Optimization and LEarning'', code: 2022ZK5ME7 MUR D.D. financing decree n. 20428 of November 6th, 2024 (CUP B53C24006410006).}

{\small\bmhead{Conflict of interest} The authors have no competing interests to declare that are relevant to the content of
this article.}

\bibliography{bibkem}%

\end{document}